\documentclass[12pt]{article}
\usepackage{amsmath,amsfonts,amsthm,mathrsfs}

\newcommand{\R}{{\mathbb R}} \newcommand{\N}{{\mathbb N}}
\newcommand{\K}{{\mathbb K}} \newcommand{\Z}{{\mathbb Z}}
  
\newcommand{\Prm}{{\mathbb P}}

\renewcommand{\epsilon}{\varepsilon } 
 
\renewcommand{\rho}{\varrho } 
\renewcommand{\phi}{\varphi }

\newcommand{\E}{{\mathbb E}\,}
\newcommand{\EE}{{\mathbb E}}

\newcommand{\ran}{{\rm ran }}
\newcommand{\ranno}{{\rm ran-non }}

\newcommand{\de}{{\rm det }}
\newcommand{\ca}{{\rm card}}

\newtheorem{theorem}{Theorem}[section]
\newtheorem{lemma}[theorem]{Lemma}
\newtheorem{corollary}[theorem]{Corollary}
\newtheorem{proposition}[theorem]{Proposition}

\begin {document}
 \title{
Randomized Complexity of Parametric Integration and the Role of Adaption II. Sobolev Spaces}

\author {Stefan Heinrich\\
Department of Computer Science\\
RPTU Kaiserslautern-Landau\\
D-67653 Kaiserslautern, Germany}  
\date{\today}%{}
\maketitle

\begin{abstract} 
We study the complexity of randomized computation of integrals depending on a parameter, with integrands from 
Sobolev spaces. That is, for $r,d_1,d_2\in\N$, $1\le p,q\le \infty$, $D_1= [0,1]^{d_1}$, and $D_2= [0,1]^{d_2}$
we are given  $f\in W_p^r(D_1\times D_2)$ and we seek to approximate
$$
Sf=\int_{D_2}f(s,t)dt\quad (s\in D_1),
$$
with error measured in the $L_q(D_1)$-norm. Our results extend previous work of Heinrich and Sindambiwe (J.\ Complexity, 15 (1999), 317--341) for $p=q=\infty$ and Wiegand (Shaker Verlag, 2006) for $1\le p=q<\infty$.
Wiegand's analysis was carried out under the assumption that $W_p^r(D_1\times D_2)$ is continuously embedded in $C(D_1\times D_2)$ (embedding condition).
We also study the case that the embedding condition
does not hold.	
For this purpose a new ingredient is developed --  a stochastic discretization technique. 

The paper is based on Part I, where vector valued mean computation -- the finite-dimensional counterpart of parametric integration -- was studied.
 In  Part I a basic problem of Information-Based Complexity on the power of adaption for linear problems in the randomized setting was solved. Here a further aspect of this problem is settled.
\end{abstract}

\section{Introduction}
\label{sec:1}
Parametric integration is the following problem. We are given a function $f$ on $D=D_1\times D_2$ and want to
compute (approximately) 
$$
Sf=\int_{D_2}f(s,t)dt\quad (s\in D_1),
$$
that is, an integral which depends on a parameter (precise definitions are given below). Here $D_1$ is the 
parameter domain and $D_2$ is the integration domain. If either of them becomes a single point, we have the integration or
the approximation problem, respectively. In general, parametric integration is an intermediate problem between 
integration and approximation,
and as such, shows features of both and clarifies the passage from one to the other. 

We consider parametric integration for $D_1=[0,1]^{d_1}$, $D_2=[0,1]^{d_2}$,  $f\in W_p^r([0,1]^{d_1+d_2})$, the Sobolev space,  and the target space (in which
the error is measured) is
$L_q([0,1]^{d_1})$. 
For the first time parametric integration was studied from the point view of information-based complexity theory in 
Heinrich, Sindambiwe \cite{HS99}, where the $p=q=\infty$ case was settled. The case $1\le p=q<\infty$ was studied by Wiegand in \cite{Wie06}.
The present paper extends the above results to the case $p\ne q$. Moreover,  \cite{Wie06} was based on the assumption that  $W_p^r([0,1]^{d_1+d_2})$ is continuously embedded into $C([0,1]^{d_1+d_2})$, the space of continuous functions (embedding condition). Our analysis also covers the case of non-embedding.

A main tool is the discretization technique for problems with standard information, which was developed in \cite{Hei03} and  \cite{Hei04b} for the analysis of integration and approximation in the quantum setting 
of information-based complexity theory. It was also applied to construct Monte Carlo methods for integration using few random bits in \cite{HNP04, GYW06, Wie06, YH08}. This method reduces an infinite dimensional problem to a series of its finite dimensional counterparts. It is deterministic and needs the embedding condition. To get away from this assumption we develop here a stochastic discretization technique which also covers the case of non-embedding. 

In Part I the adaption problem in the randomized setting was settled: Is there a constant $c>0$ such that for all linear problems
$\mathcal{P}=(F,G,S,K,\Lambda)$ and all $n\in\N$
\begin{equation*}
e_n^\ranno (S,F,G)\le ce_n^\ran (S,F,G) \, {\bf ?}
\end{equation*}
see \cite{Nov96}, \cite{NW08}, and Part I for further references.
It was shown in Part I that there is a sequence of vector-valued mean computation problems $\mathcal{P}_n$ $(n\in\N)$  such that
\begin{equation}
\label{AK0}
\lim_{n\to \infty}\frac{e_n^\ranno (S_n,F_n,G_n)}{e_n^\ran (S_n,F_n,G_n)} = \infty.
\end{equation}
The results of the present paper imply that there is a single problem $\mathcal{P}$ such that 
\begin{equation}
\label{AK1}
\lim_{n\to \infty}\frac{e_n^\ranno (S,F,G)}{e_n^\ran (S,F,G)} = \infty.
\end{equation}
We show that the parametric integration problem with $2<p<q\le \infty$ has this property. Moreover, there are instances of this problem where the quotient in \eqref{AK1}, that is, the gap between non-adaptive and adaptive randomized minimal errors is (up to log factors) of order $n^{1/8}$. For problems with larger gaps (in the sense of \eqref{AK0}) we refer to \cite{Hei23b}.

The paper is organized as follows: In Section \ref{sec:2} we present the needed prerequisites including some results from Part I. In Section \ref{sec:3} we develop the randomized discretization technique (which contains the deterministic one as a special case). This technique is used in Section \ref{sec:4} to set up and analyze non-adaptive and adaptive algorithms for parametric integration with error of optimal order (often up to logarithms). Lower bounds and the main theorem on the complexity are proved in Section \ref{sec:5}. The final Section \ref{sec:6} contains a discussion of the maximal gap between non-adaptive and adaptive randomized algorithms.

\section{Preliminaries}
\label{sec:2}
 We denote ${\mathbb{N}}=\{1,2,\dots\}$,  ${\mathbb{N}}_0={\mathbb{N}}\cup\{0\}$, and for $N\in\N$,   $\Z[1,N]:=\{1,2,\dots,N\}$. The symbol ${\mathbb{K}}$ stands for the scalar field ${\mathbb R}$ or ${\mathbb C}$.
We often use the same symbol
$c, c_0,c_1,\dots$ for possibly different constants, even if they appear in a sequence
of relations. However, some constants are supposed to have the same meaning throughout the paper  --  these are denoted by symbols $c(0),c(1),\dots$. Throughout this paper $\log$ means $\log_2$. The unit ball of a normed or semi-normed space $X$ is denoted by $B_X$, the $\sigma$-algebra of Borel subsets of a normed space $X$ by $\mathcal{B}(X)$ and the space of bounded linear operators from $X$ to a normed space $Y$ by $L(X,Y)$.
 
Let $d\in \N$, $D=[0,1]^d$.  By $C(D)$ we denote the space of continuous functions on $D$,
equipped with the supremum norm. For $1\leq p \le \infty$ let $L_p(D)$ be the
space of equivalence classes of $\K$-valued Borel measurable functions which are $p$-integrable with respect to the Lebesgue measure, endowed
with the usual norm
\[
\| f \|_{L_p(D)} = \left\{\begin{array}{lll}
   &  \left( \int_D |f(t)|^p dt \right)^{1/p}\quad\mbox{if}\quad  p<\infty  \\[.2cm]
   &  \mbox{\rm ess-sup}_{t\in D} |f(t)|\quad\mbox{if}\quad p=\infty.
    \end{array}
\right.   
\]
For $r\in\N$ the Sobolev space $W_p^r(D)$ consists of all  $f\in L_p(D)$ such that
for all\\
 $\alpha =(\alpha_1, \dots , \alpha_d) \in \N_0^d$ with
$|\alpha|:= \sum_{j=1}^d \alpha_j \leq r$, the generalized partial derivative
$\partial^{\alpha} f$ belongs to $L_p(D)$. The norm on $W_p^r(D)$ is defined by
\[
\| f \|_{W_p^r(D)} = \left\{\begin{array}{lll}
   & \left( \sum_{|\alpha|\leq r} \| \partial^{\alpha} f \|_{L_p(D)}^p \right)^{1/p}\quad\mbox{if}\quad p<\infty    \\
   & \max_{|\alpha|\leq r} \| \partial^{\alpha} f \|_{L_\infty(D) }\quad \mbox{if}\quad p=\infty.
    \end{array}
\right.
\]
For $r=0$ we set $W_p^r(D)=L_p(D)$. 
We  recall from  \cite{Ada75}, Th.\ 5.4,  that $W_p^r(D)$ is continuously embedded into  $C(D)$ if and only if
\begin{equation}\label{PA2}
\left.\begin{array}{lll}
&p=1 &  \;\mbox{and} \quad r/d\ge 1   \\
\mbox{or}\\
  &1<p\le\infty & \;\mbox{and} \quad  r/d>1/p
    \end{array}\right\} 
\end{equation}
(embedding condition). 

Now let $1\le p,q\le \infty$,  $d_1, d_2 \in \N_0$, with $d=d_1+d_2$. First assume $d_1,d_2\in\N$ and let $D_1=[0,1]^{d_1}, D_2=[0,1]^{d_2}$, thus
$D= D_1 \times D_2$.
Consider the operator of parametric integration 
$S: W_p^r(D)\to  L_q(D_1)$ defined by
\begin{equation}
\label{AI3}
(Sf)(s) = \int_{D_2} f(s,t)\, dt \quad(s\in D_1).
\end{equation}
We always assume $d\ge 1$, but we also want to include the border-line cases $d_1=0$ and $d_2=0$.
In the case of $d_\iota=0$ $(\iota=1,2)$ we let  $D_\iota=\{0\}$ be the one-point-set. The integral over $D_\iota$ as well as the space $L_q(D_\iota)$ are considered as defined  with respect to  the trivial probability measure on  $D_\iota=\{0\}$. In other words,  for $d_1=0$ we have 
$L_q(D_1)=\K$ and $S: W_p^r(D)\to  \K$, 
\begin{eqnarray*}
 Sf = \int_D f(t)\, dt 
\end{eqnarray*}
is the integration operator. It is well-defined and continuous for all $r\in\N_0$ and $1\le p\le \infty$.

If $d_2=0$ we have  
$\int_{D_2} f(s)\, dt=f(s)$ and therefore $S: W_p^r(D)\to  L_q(D)$, 
\begin{eqnarray}
\label{AI5}
 (Sf)(s) =f(s) \quad(s\in D)
\end{eqnarray}
is the approximation operator. 
If $d_1\ge 1$, the operator $S$ from \eqref{AI3} and \eqref{AI5} is well-defined and continuous iff the embedding of $W_p^r(D_1)$ to  $L_q(D_1)$ is continuous, that is, iff
\begin{equation}\label{B1}
\left.\begin{array}{llll}
&1\le q<\infty & &  \;\mbox{and} \quad \frac{r}{d_1}\ge \left(\frac{1}{p}-\frac{1}{q}\right)_+  \\
\mbox{or}\\
&q=\infty,  &1<p<\infty, & \;\mbox{and} \quad  \frac{r}{d_1}>\frac{1}{p}\quad  
   \\\mbox{or}\\
&q=\infty, &p\in\{1,\infty\}, &  \;\mbox{and} \quad\frac{r}{d_1}\ge\frac{1}{p}. 
    \end{array}\right\} 
\end{equation}
Here we used the notation $a_+=\max(a,0)$ for $a\in\R$.

In the sequel we also need the stronger condition (of compact embedding) 
\begin{equation}\label{BB1}
 \frac{r}{d_1}>\left(\frac{1}{p}-\frac{1}{q}\right)_+ . 
\end{equation}

In the terminology of Section 2 of Part I we put $G=L_q(D_1)$, and $K=\K$.  
We consider standard information, that is, function values. This needs some care since $W_p^r(D)$ consists of equivalence classes of functions.

If \eqref{PA2} holds, then
each class $[f]\in W_p^r(D)$ contains a unique continuous representative and we set
$\mathcal{W}_p^r(D)=\{f\in C(D): [f]\in W_p^r(D)\}$, equipped with the norm $\|f\|_{\mathcal{W}_p^r(D)}:=\|[f]\|_{W_p^r(D)}$. 

If \eqref{PA2}  does not hold, we let $\mathcal{W}_p^r(D)$ be the respective Sobolev space of functions (not equivalence classes), thus $f\in \mathcal{W}_p^r(D)$ iff $[f]\in W_p^r(D)$, where $[f]$ is the equivalence class of $f$ with respect to equality up to a subset of $D$ of Lebesgue measure zero. This is a linear space and   $\|f\|_{\mathcal{W}_p^r(D)}:=\|[f]\|_{W_p^r(D)}$ is a semi-norm on it. Clearly, $S$ can also be viewed as an operator from $\mathcal{W}_p^r(D)$ to $L_q(D_1)$.

Finally, in both cases we put 
$F=B_{\mathcal{W}_p^r(D)}$ or $F=\mathcal{W}_p^r(D)$ and $\Lambda=\{\delta_t:\,t\in D\}$, where  $\delta_t(f)=f(t)$.  Consequently, here we study the problems
\begin{eqnarray}
&&(B_{\mathcal{W}_p^r(D)},L_q(D_1), S, \Lambda, \K)\label{AL5}
\\
&&(\mathcal{W}_p^r(D),L_q(D_1), S, \Lambda, \K).\label{AL6}
\end{eqnarray}

We also recall the discrete version of parametric integration studied in Part I.
Let $M,M_1,M_2$ be finite sets and define the space
$L_p^M$ to be the set of all functions $f: M \rightarrow \K$
with norm
\[
\| f \|_{L_p^M} = \left\{\begin{array}{lll}
 \displaystyle \left( \frac{1}{|M|} \sum_{i\in M} |f(i)|^p \right)^{1/p} & \quad\mbox{if}\quad  p<\infty  \\[.5cm]
  \displaystyle\max_{i\in M} |f(i)| & \quad\mbox{if}\quad p=\infty
    \end{array}
\right.
\]
and the operator of discrete parametric integration (or vector-valued mean computation) as 
\begin{equation}
\label{eq:3}
S^{M_1,M_2} : L_p^{M_1\times M_2} \rightarrow L_q^{M_1}
\end{equation}
with
\begin{equation}
\label{AI1}
(S^{M_1,M_2}f)(i)  = \frac{1}{|M_2|} 
\sum_{j\in M_2} f(i,j) .
\end{equation}
In this connection we consider the problems
\begin{eqnarray}
&&\left( B_{L_p^{M_1\times M_2}},L_q^{M_1},S^{M_1,M_2},\K,\Lambda\right)\label{AL7}
\\
&&\left( L_p^{M_1\times M_2},L_q^{M_1},S^{M_1,M_2},\K,\Lambda\right),
\label{AL8}
\end{eqnarray}
where
$\Lambda=\{\delta_{ij}:\, i\in M_1,\,j\in M_2\}$ with $\delta_{ij}(f)=f(i,j)$.
Given $N,N_1,N_2\in\N$, we  write $L_p^{N}$ for $L_p^{\Z[1,N]}$, 
$L_p^{N_1,N_2}$ for $L_p^{\Z[1,N_1]\times \Z[1,N_2]}$, and $S^{N_1,N_2}$ for $S^{\Z[1,N_1],\Z[1,N_2]}$.
Furthermore, for  $1\le p,q\le \infty$  we define 
\begin{eqnarray}
\bar{p}=\min(p,2),\quad
\sigma_1=\left\{\begin{array}{lll}
	1 & \mbox{if} \quad p=q=\infty   \\
	0 & \mbox{otherwise.}    \\
	\end{array}
	\right.
\label{Q8}
\end{eqnarray}
For the definition of the non-adaptive randomized algorithm   $A_{n}^{(2)}=\big(A_{n,\omega}^{(2)}\big)_{\omega\in\Omega}$ and the adaptive randomized algorithm $A_{n,m}^{(3)}=\big(A_{n,m,\omega}^{(3)}\big)_{\omega\in\Omega}$  for problem \eqref{AL8} we refer to Part I. The next two results are Proposition 4.2 and 4.3 of Part I. Since they will be used in connection with \eqref{WM6} below for several $n$ simultaneously, we denote the involved probability space by $(\Omega_2,\Sigma_2,\Prm_2)$ and assume w.l.o.g. that it is the same for all $n\in\N$ (the index 2 is convenient for the further notation).  The respective expectation is denoted by $\EE_2$.
\begin{proposition}
	\label{pro:4}
	Let $1 \le p,q \le \infty $ and $1\le w\le p$, $w<\infty$. Then there is a constant $c>0$ such that for all $n,N_1,N_2\in \N$ with $n<N_1N_2$ and all $f\in L_p^{N_1,N_2}$ 
	\begin{equation}
	\label{K1}
	\EE_2 A_{n,\omega_2}^{(2)}f= S^{N_1,N_2}f\quad(n\ge N_1),\quad \ca(A_{n,\omega_2}^{(2)})\le 2n\quad(\omega_2\in\Omega_2),
	\end{equation}
and
\begin{eqnarray}
	\label{K2}
		\lefteqn{
		 \left(\EE_2 \|S^{N_1,N_2}f-A_{n,\omega_2}^{(2)}f\|_{L_q^{N_1}}^w\right)^{1/w}}
		\nonumber\\[.2cm]
		&\le &
	cN_1^{\left(1/p-1/q\right)_+}
		\left\lceil\frac{n}{N_1}\right\rceil^{-1+1/\bar{p}}
		\min\left(\log (N_1+1),\left\lceil\frac{n}{N_1}\right\rceil\right)^{\sigma_1/2}\|f\|_{L_p^{N_1,N_2}}.
	\end{eqnarray}
\end{proposition}
The adaptive algorithm $A_{n,m}^{(3)}$ will only be used in the case $2<p<q\le \infty$. The following is  Proposition 4.3 of Part I.
\begin{proposition}
\label{pro:7}
Let $2< p <q\le  \infty$ and $1\le w<\infty$. Then there exist constants $c_1,c_2>0$ such that the following hold for all $m,n,N_1,N_2\in\N$ and $ f\in L_p^{N_1,N_2}$:
\begin{equation}
\label{WM1}
\ca(A_{n,m,\omega_2}^{(3)})\le 6mn
\end{equation}
and for $m\ge c_1\log(N_1+N_2)$, 
$1 \le n < N_1N_2$ 
\begin{eqnarray}
\label{mc-eq:2}
\lefteqn{\left(\EE_2 \|S^{N_1,N_2}f-A_{n,m,\omega_2}^{(3)}f\|_{L_q^{N_1}}^w\right)^{1/w}}
\notag\\
&\le& 
 c_2\Bigg(N_1^{1/p-1/q}\left\lceil\frac{n}{N_1}\right\rceil^{-\left(1-1/p\right)}+\left\lceil\frac{n}{N_1}\right\rceil^{-1/2}\Bigg)\|f\|_{L_p^{N_1,N_2}}.\quad
\end{eqnarray}
\end{proposition}
Let us also note that for all $\omega_2\in\Omega_2$,
\begin{eqnarray}
A_{n,\omega_2}^{(2)}(0)=0, \quad A_{n,m,\omega_2}^{(3)}(0)=0.
\label{AG7}
\end{eqnarray}
This follows directly from the definitions of $A_{n,\omega_2}^{(2)}$ and  $A_{n,m,\omega_2}^{(3)}$ in Part I.
The next proposition summarizes the lower bounds in Theorem 4.5 of Part I.
\begin{proposition}
\label{pro:1}
Let $1 \leq p,q \le \infty $. Then there exists constants $0<c_0<1$, $c_1,\dots,c_4>0$,  such that for $n,N_1,N_2\in\N$  with 
$n < c_0N_1N_2$ the following hold: \\
If $p\le 2$ or $p\ge q$, then  
\begin{eqnarray}
\label{AD2}
\lefteqn{e_n^{\rm ran-non}(S^{N_1,N_2},B_{L_p^{N_1,N_2}},L_q^{N_1})\ge e_n^\ran(S^{N_1,N_2},B_{L_p^{N_1,N_2}},L_q^{N_1})}
\nonumber\\[.2cm]
&\ge&
c_1N_1^{\left(1/p-1/q\right)_+}\left\lceil\frac{n}{N_1}\right\rceil^{-\left(1-1/\bar{p}\right)}\left( \min \left( \log (N_1+1), \left\lceil\frac{n}{N_1}\right\rceil \right) \right)^{\sigma_1/2}.
\end{eqnarray}
If $2<p<q$, then  
\begin{eqnarray}
\label{AD3}
\lefteqn{e_n^\ran(S^{N_1,N_2},B_{L_p^{N_1,N_2}},L_q^{N_1})}
\nonumber\\ 
&\ge&c_2N_1^{1/p-1/q}\left\lceil\frac{n}{N_1}\right\rceil^{-\left(1-1/p\right)}+c_2\left\lceil\frac{n}{N_1}\right\rceil^{-1/2}  \left( \log (N_1+1) \right)^{\delta_{q,\infty}/2}
\end{eqnarray}
and
\begin{eqnarray}
e_n^{\rm ran-non}(S^{N_1,N_2},B_{L_p^{N_1,N_2}},L_q^{N_1})\ge c_3N_1^{1/p-1/q}\left\lceil\frac{n}{N_1}\right\rceil^{-1/2}
 .
\label{AD4}
\end{eqnarray}
In the deterministic setting we have 
\begin{eqnarray}
\label{AD5}
e_n^\de(S^{N_1,N_2},B_{L_p^{N_1,N_2}},L_q^{N_1})\ge c_4N_1^{\left(1/p-1/q\right)_+}.
\end{eqnarray}
\end{proposition}

\section{Deterministic and stochastic discretization}
\label{sec:3}
In this section we are concerned with discretization, that is, reducing the problem of parametric integration to a family of vector-valued summation problems. 
We present the needed ingredients of the discretization technique developed in  \cite{Hei03, Hei04b} and used in \cite{HNP04,GYW06, Wie06,YH08}. This technique is based on deterministic function evaluations, the availability of which requires that $W_p^r(D)$ is continuously embedded in $C(D)$. In this paper we also study the case of non-embedding. For this purpose we develop here, in addition, a stochastic discretization technique. The presentation below also includes deterministic discretization, which appears as a special case of the stochastic one. In this section we assume $d_1,d_2\in\N_0$, $d_1+d_2\ne 0$, that is, we include the border-line cases $d_1=0$, $d_2\ne 0$ (integration) and $d_1\ne 0$, $d_2=0$ (approximation).

For $l\in \N_0$ let 
\begin{equation*}
D = \bigcup_{i=1}^{2^{dl} }D_{li}
\end{equation*}
be a partition into congruent
cubes of disjoint interior with side length $2^{-l}$. If $d_1d_2\ne 0$, we assume in addition, that 
\begin{equation*}
 D_{li}=D_{li_1}^{(1)}\times D_{li_2}^{(2)}\quad i=2^{d_2l}(i_1-1)+i_2,
\end{equation*}
where
\[
D_1  = \bigcup_{i_1=1}^{2^{d_1 l}} D^{(1)}_{li_1}, \quad D_2  = \bigcup_{i_2=1}^{2^{d_2 l}} D^{(2)}_{li_2}
\]
are respective partitions of $D_1$ and $D_2$.

Let $s_{li}\in D_{li}$ denote the point with minimal Euclidean norm.
We introduce the following operators 
$E_{li}$ and $R_{li}$ on $\mathcal{F}(D,\K)$, the set
of all $\K$-valued functions on $D$,  by setting for 
$s \in D$
\[
(E_{li} f)(s) = f(s_{li} +2^{-l}s) 
\]
and
\begin{equation}
\label{WH4}
(R_{li} f)(s) = \left\{ \begin{array}{ll}
f(2^l(s-s_{li})) & \text{ if } s\in D_{li}\\
0 & \text{ otherwise.} 
\end{array}
\right.
\end{equation}
We also need the operators $R^{(1)}_{li_1}$ for $1\le i_1\le 2^{d_1 l}$, 
which are defined analogously, with $D_1$ instead of $D$. Definition \eqref{WH4} also makes sense if $d_1=0$. Then we have $\mathcal{F}(\{0\},\K)=\K$ and $R^{(1)}_{li_1}=Id_\K$, where $Id$ denotes the identity operator.

Fix $r\in\N$ and $0<\delta<1$. If $r\ge 2$,  let $(t_k)_{k=1}^{\kappa}$ be the uniform grid of mesh size $(r-1)^{-1}(1-\delta)$ on $[0,1-\delta]^d$ and if $r=1$, set $\kappa=1$ and $t_1=0$. Let $P$ be for $d=1$ the respective Lagrange interpolation operator of degree $r-1$ and for $d>1$ its tensor product. We consider $P$ as an operator from $\mathcal{F}(D,\K)$ to $L_\infty(D)$. 
Let $\mathcal{P}^{r-1}_{\rm max}(Q)$ denote the space of polynomials of maximum degree $r-1$ on a set $Q\subset \R^d$.
Then $P$ can be represented as
\begin{equation}
\label{AF8}
(Pf)(t) = \sum_{k=1}^{\kappa} f(t_k) \varphi_k \quad (f\in \mathcal{F}(D,\K), t\in D).
\end{equation}
with the $d$-dimensional tensor product Lagrange polynomials $(\phi_k)_{k=1}^{\kappa}$, which form a basis in the space $\mathcal{P}^{r-1}_{\rm max}(\R^d)$. 

Now we randomize this operator. Let $\rho$ be a random variable over a probability space $(\Omega_1,\Sigma_1,\Prm_1)$ with values in $[0,1]^d$. 
Define 
\begin{equation}
(P_{0,\rho}f)(t) = \chi_{D}(t)\sum_{k=1}^{\kappa} f(t_k+\delta\rho) \varphi_k(t-\delta\rho) \quad (f\in \mathcal{F}(D,\K), t\in D).
\label{W6}
\end{equation}
Expanding $\varphi_k(t-\delta\rho)$ with respect to $t$ and $\rho$, we can represent 
$$
\varphi_k(t-\delta\rho)=\sum_{j=1}^{\kappa}a_{jk}(\rho)\phi_j(t),
$$
with $a_{jk}\in \mathcal{P}^{r-1}_{\rm max}(\R^d)$, hence
\begin{equation*}
P_{0,\rho}f = \chi_{D}\sum_{j=1}^{\kappa}\sum_{k=1}^{\kappa} a_{jk}(\rho)f(t_k+\delta\rho) \varphi_j \quad (f\in \mathcal{F}(D,\K)).
\end{equation*}
For $f\in \mathcal{F}(D,\K)$ and $l\in \N_0$ we set
\begin{equation}
\label{WF1}
P_{l\rho} f = \sum_{i=1}^{2^{dl}} R_{li} P_{0,\rho} E_{li} f
= \sum_{i=1}^{2^{dl}}\sum_{j=1}^{\kappa} \sum_{k=1}^{\kappa}a_{jk}(\rho)f(s_{li}+2^{-l}(t_k+\delta\rho)) R_{li}\varphi_j.
\end{equation}
Next we define
\begin{eqnarray}
P_{0,\rho}'f  =  (P_{1,\rho} - P_{0,\rho})f 
&=&\sum_{i_0=1}^{2^{d}}\sum_{j_0=1}^{\kappa}\sum_{k_0=1}^{\kappa}a_{j_0k_0}(\rho)f(s_{1,i_0}+2^{-1}(t_{k_0}+\delta\rho)) R_{1,i_0} \varphi_{j_0}
\nonumber\\
&&-\chi_{D}\sum_{m=1}^{\kappa} \sum_{k_0=1}^{\kappa}a_{mk_0}(\rho)f(t_{k_0}+\delta\rho) \varphi_m .\label{P4}
\end{eqnarray}
Expanding 
\begin{equation*}
\chi_{D_{1,i_0}}\phi_m=\sum_{j_0=1}^{\kappa}\beta_{i_0mj_0}R_{1,i_0}\phi_{j_0}
\end{equation*}
with $\beta_{i_0mj_0}\in\K$, \eqref{P4} turns into 
\begin{eqnarray*}
P_{0,\rho}'f 
&=&\sum_{i_0=1}^{2^{d}}\sum_{j_0=1}^{\kappa}\sum_{k_0=1}^{\kappa}\Bigg(a_{j_0k_0}(\rho)f(s_{1,i_0}+2^{-1}(t_{k_0}+\delta\rho)) 
\nonumber\\
&&\hspace{2.7cm}-\sum_{m=1}^{\kappa} \beta_{i_0mj_0}a_{mk_0}(\rho)f(t_{k_0}+\delta\rho) \Bigg)R_{1,i_0} \varphi_{j_0},
\end{eqnarray*}
which we write as
\begin{equation*}
P_{0,\rho}'f = \sum_{j=1}^{\kappa'} \sum_{k=1}^{\kappa''} b_{jk}(\rho) f(t_{jk}(\rho)) 
 \psi_{j},
\end{equation*}
with 
\begin{equation}
\label{eq:24}
\kappa'=2^{d}\kappa,\quad\kappa''=2\kappa,
\end{equation}
and for $1\le i_0\le 2^d$,  $1\le j_0\le \kappa$, $j=\kappa (i_0-1)+j_0$, $1\le k\le 2\kappa$
\begin{eqnarray}
\psi_j&=&R_{1,i_0} \varphi_{j_0} 
\label{AG0}\\
t_{jk}(\rho)&=&\left\{\begin{array}{lll}
   s_{1,i_0}+2^{-1}(t_k+\delta\rho)& \quad\mbox{if}\quad  1\le k\le \kappa  \\
   t_{k-\kappa}+\delta\rho& \quad\mbox{if}\quad \kappa+1\le k\le 2\kappa
    \end{array}
\right. 
\nonumber\\[.2cm]
b_{jk}(\rho)&=&\left\{\begin{array}{lll}
   a_{j_0k}(\rho)& \quad\mbox{if}\quad  1\le k\le \kappa  \\[.2cm]
   \displaystyle\sum_{m=1}^{\kappa} \beta_{i_0mj_0}a_{m,k-\kappa}(\rho)& \quad\mbox{if}\quad \kappa+1\le k\le 2\kappa.
    \end{array}
\right. 
\nonumber
\end{eqnarray}
Thus the $t_{jk}(\rho)$ are $D$-valued and the  $b_{jk}(\rho)$  scalar-valued random variables. 

Now we put
\begin{equation}
\label{eq:25}
\psi_{lij} = R_{li} \psi_j\quad (1\le i\le 2^{dl},\,1\le j\le \kappa')
\end{equation}
and let
\begin{equation}
\label{eq:29}
P_{l\rho}'f  =  \sum_{i=1}^{2^{dl}} R_{li}P_{0,\rho}'E_{li}f =\sum_{i=1}^{2^{dl}}\sum_{j=1}^{\kappa'}
\sum_{{k}=1}^{\kappa''} b_{jk}(\rho)f(s_{li}+2^{-l}t_{jk}(\rho))\psi_{lij}.
\end{equation}
Observe that 
\begin{eqnarray}
P_{l\rho}'  =  \sum_{i=1}^{2^{dl}} R_{li}\Bigg(\sum_{i_0=1}^{2^d}R_{1,i_0}P_{0,\rho}E_{1,i_0}-P_{0,\rho}\Bigg)E_{li}
=P_{l+1,\rho}-P_{l\rho},
\label{W8}
\end{eqnarray}
hence for $l_0,l_1 \in \N_0$, $l_0\le l_1$ 
\begin{equation}
\label{C1}
P_{l_1\rho}=P_{l_0\rho}+\sum_{l=l_0}^{l_1-1}P'_{l\rho}.
\end{equation}
For the sake of brevity we denote 
\begin{equation*}
\mathcal{I}_{1,l}=\Z[1,2^{d_1l}]\times\Z[1,\kappa'], \quad \mathcal{I}_{2,l}=\Z[1,2^{d_2l}]
\end{equation*}
and set
\begin{equation}
\label{eq:26}
N_{1,l}:=|\mathcal{I}_{1,l}|=\kappa' 2^{d_1l},\quad N_{2,l}:=|\mathcal{I}_{2,l}|=2^{d_2l}. 
\end{equation}
Define the operator 
$
T_l:{\rm span}\left\{\psi_{lij}:i\in\Z[1,2^{dl}],
j\in\Z[1,\kappa']\right\}\to L_p^{\mathcal{I}_{1,l}\times \mathcal{I}_{2,l}}$ by  
$$
  T_l\psi_{lij} = \delta_{(i_1,j),i_2}\quad((i_1,j)\in \mathcal{I}_{1,l},\, i_2\in \mathcal{I}_{2,l},\, i=2^{d_2l}(i_1-1)+i_2).
$$
It follows from the definition \eqref{AG0}, \eqref{eq:25} of $\psi_{lij}$ and the independence of the sequence $(\phi_k)_{k=1}^\kappa$, see \eqref{AF8}, that $T_l$ is correctly defined and 
\begin{equation}
\|T_lf\|_{L_p^{\mathcal{I}_{1,l}\times \mathcal{I}_{2,l}}}\le c\|f\|_{L_p(D)}\qquad \left(f\in {\rm span}\left\{\psi_{lij}\right\}\right). 
\label{S3}
\end{equation}
We define furthermore
\begin{equation}
\label{X1}
U_{l\rho}: \mathcal{W}_p^r(D) \rightarrow L_p^{\mathcal{I}_{1,l}\times \mathcal{I}_{2,l}}, \quad U_{l\rho}=T_lP_{l\rho}'.
\end{equation}
Then
\begin{equation}
\label{WF3}
(U_{l\rho}f)((i_1,j),i_2)=\sum_{{k}=1}^{\kappa''} b_{jk}(\rho) f(s_{li} +2^{-l}t_{jk}(\rho)).
\end{equation}

Next we turn to the operator of parametric integration.   Put
$$
\theta_j=S\psi_j\in L_\infty(D_1)\quad(1\le j\le\kappa')
$$
and 
\begin{equation}
\label{WF6}
\theta_{li_1j}=R^{(1)}_{li_1}\theta_j\quad ( (i_1,j)\in\mathcal{I}_{1,l}).
\end{equation}
Observe that for $(i_1,j)\in\mathcal{I}_{1,l}, \,i_2\in \mathcal{I}_{2,l},\, i=2^{d_2l}(i_1-1)+i_2$,
\begin{equation}
\label{P2}
S\psi_{lij}=SR_{li} \psi_j=2^{-d_2l} R_{li_1}^{(1)}S\psi_j=2^{-d_2l}\theta_{li_1j}.
\end{equation}
Let 
\begin{equation}
\label{eq:37}
V_l:  L_q(\mathcal{I}_{1,l}) \rightarrow L_q(D_1),\quad V_l g = \sum_{i_1=1}^{2^{d_1l}} \sum_{j=1}^{\kappa'}
g(i_1,j) \theta_{li_1j} .
\end{equation}
According to \eqref{WF6}, for $d_1\ge 1$ the supports of the $\theta_{li_1j}$ are disjoint for different $i_1$, therefore 
\begin{eqnarray}
\|V_lg\|_{L_q(D_1)}&=&\Bigg(\sum_{i_1=1}^{2^{d_1l}} \bigg\|\sum_{j=1}^{\kappa'}g(i_1,j) \theta_{li_1j} \bigg\|_{L_q(D_1)}^q\Bigg)^{1/q}
\nonumber\\
&=&\Bigg(2^{-d_1l}\sum_{i_1=1}^{2^{d_1l}} \bigg\|\sum_{j=1}^{\kappa'}g(i_1,j) \theta_j \bigg\|_{L_q(D_1)}^q\Bigg)^{1/q}
\le c \|g\|_{L_q(\mathcal{I}_{1,l})},
\label{P8}
\end{eqnarray}
with the obvious modifications for $q=\infty$. This relation trivially also holds for $d_1=0$, since in this case by \eqref{WF6} $\theta_{li_1j}=\theta_j$. Now  \eqref{eq:29}, \eqref{WF3}, \eqref{P2},  and \eqref{eq:37} give
\begin{eqnarray}
S P_{l\rho}'f  &=&\sum_{i=1}^{2^{dl}}\sum_{j=1}^{\kappa'}
\sum_{{k}=1}^{\kappa''} b_{jk}(\rho)f(s_{li}+2^{-l}t_{jk}(\rho))S\psi_{lij}
\notag\\
&=&\sum_{i_1=1}^{2^{d_1l}}\sum_{j=1}^{\kappa'}2^{-d_2l}\sum_{i_2=1}^{2^{d_2l}}
\sum_{{k}=1}^{\kappa''} b_{jk}(\rho)f(s_{li}+2^{-l}t_{jk}(\rho))\theta_{li_1j}
\notag\\
&=&\sum_{i_1=1}^{2^{d_1l}}\sum_{j=1}^{\kappa'}2^{-d_2l}\sum_{i_2=1}^{2^{d_2l}}
(U_{l\rho}f)((i_1,j),i_2)\theta_{li_1j}
\notag\\
&=&\sum_{i_1=1}^{2^{d_1l}}\sum_{j=1}^{\kappa'}(S^{\mathcal{I}_{1,l}, \mathcal{I}_{2,l}}
U_{l\rho}f)(i_1,j)\theta_{li_1j}=V_lS^{\mathcal{I}_{1,l}, \mathcal{I}_{2,l}}U_{l\rho}f,
\label{AG4}
\end{eqnarray}
where $S^{\mathcal{I}_{1,l}, \mathcal{I}_{2,l}}$ is the operator of  vector valued mean computation defined in \eqref{eq:3} and \eqref{AI1}.
Using \eqref{C1}, this implies
\begin{equation}
\label{eq:41}
S = S-SP_{l_1\rho} + SP_{l_0\rho} +\sum_{l=l_0}^{l_1-1} SP_{l\rho}'
=S-SP_{l_1\rho}+SP_{l_0\rho} +\sum_{l=l_0}^{l_1-1} V_l S^{\mathcal{I}_{1,l}, \mathcal{I}_{2,l}} U_{l\rho}.
\end{equation}
This way 
the problem $S$ is reduced 
to the computation of $SP_{l_0\rho}$, which will be done exactly, provided $l_0$ is suitable chosen, and to the (approximate) computation of the vector-valued means $S^{\mathcal{I}_{1,l}, \mathcal{I}_{2,l}}$ for $l_0\le l<l_1$. 
So let $l_0\le l_1\in \N_0$ and let  $A_{l,\omega}:L_p^{\mathcal{I}_{1,l}, \mathcal{I}_{2,l}}\to L_q^{\mathcal{I}_{1,l}}$ for each $l$ with $l_0\le l<l_1$ be a randomized algorithm for problem  \eqref{AL8}, independent of $\rho$. For convenience of presentation, we assume that the  $A_{l,\omega}$ are defined over a different from $(\Omega_1,\Sigma_1,\Prm_1)$ probability space $(\Omega_2,\Sigma_2,\Prm_2)$ and put 
\begin{equation*}
(\Omega,\Sigma,\Prm)=(\Omega_1,\Sigma_1,\Prm_1)\times(\Omega_2,\Sigma_2,\Prm_2).
\end{equation*}
The expectation with respect to $\Prm_2$ and $\Prm$ are denoted by $\EE_2$ and $\EE$, respectively. 
Now we define a randomized  algorithm $A=(A_\omega)_{\omega\in\Omega}$ for problem \eqref{AL6} by  setting 
\begin{equation}
\label{WM6}
A_\omega(f)=SP_{l_0,\rho(\omega_1)}f +\sum_{l=l_0}^{l_1-1} V_lA_{l,\omega_2} ( U_{l,\rho(\omega_1)} f)\quad(\omega=(\omega_1,\omega_2)).
\end{equation}
For fixed $\omega\in\Omega$,  $A_\omega$ is indeed  a deterministic algorithm in the formal sense of Section 2 of Part I. This follows from  Lemma 2 and 3 of \cite{Hei05b}. Moreover, if the $A_{l,\omega_2} $ are non-adaptive, then so is $A_\omega$. Finally, we assume that for $l_0\le l<l_1$  and  $f\in \mathcal{W}_p^r(D)$ the mappings
\begin{eqnarray}
(\omega_1,\omega_2)\to\left\{\begin{array}{lll}
   \ca\big(A_{l,\omega_2},U_{l,\rho(\omega_1)} f\big)   \\[,2cm]
   A_{l,\omega_2} ( U_{l,\rho(\omega_1)} f)   
    \end{array}
\right. \label{AM5}
\end{eqnarray}
are $\Sigma$-measurable. Then $(A_\omega)_{\omega\in\Omega}$ is a randomized algorithm for problem \eqref{AL6}.
\begin{proposition}
\label{pro:5}
Let $r\in \N_0$, $1 \le p,q \le \infty $, $1\le w<\infty$, and assume that \eqref{B1} and \eqref{AM5} hold.  Then
there is a constant $c>0$ such that  for all $l_0,l_1$, $\rho$ as above, and $f\in \mathcal{W}_p^r(D)$
\begin{eqnarray}
	\label{WO1}
			\lefteqn{\left(\E \|Sf-A_{\omega}(f)\|_{L_q(D_1)}^w\right)^{1/w}}
\nonumber\\
		&\le& c(\EE_1\|Sf-SP_{l_1\rho(\omega_1)}f\|_{L_q(D_1)}^w)^{1/w}+ c\sum_{l=l_0}^{l_1-1}\left(\EE_2\|A_{l,\omega_2}(0)\|_{L_q^{\mathcal{I}_{1,l}}}^w\right)^{1/w}
\nonumber\\
&&+ c\sum_{l=l_0}^{l_1-1}\sup_{g\in L_p^{\mathcal{I}_{1,l}\times \mathcal{I}_{2,l}}\setminus\{0\}}\left(\|g\|_{L_p^{\mathcal{I}_{1,l}\times \mathcal{I}_{2,l}}}^{-w}\EE_2\|S^{\mathcal{I}_{1,l}, \mathcal{I}_{2,l}}g-A_{l,\omega_2}(g)\|_{L_q^{\mathcal{I}_{1,l}}}^w\right)^{1/w}
\nonumber\\
&&\quad\times\left(\EE_1\left\|(P_{l+1,\rho(\omega_1)}-P_{l\rho(\omega_1)})f\right\|_{L_p(D)}^w\right)^{1/w}.
\end{eqnarray}
Moreover, 
\begin{equation}
\label{U4}
\ca(A_\omega,\mathcal{W}_p^r(D))\le \kappa 2^{(d_1+d_2)l_0}+\kappa''\sum_{l=l_0}^{l_1-1} \ca\big(A_{l,\omega_2},L_p^{\mathcal{I}_{1,l}\times \mathcal{I}_{2,l}}\big)\quad (\omega=(\omega_1,\omega_2)\in\Omega).
\end{equation}
\end{proposition}
\begin{proof} Relation \eqref{U4} follow directly from \eqref{WF1} and \eqref{WF3}. 
By \eqref{eq:41} and \eqref{WM6} 
\begin{eqnarray}
Sf-A_{\omega}(f)
&=&Sf-SP_{l_1\rho}f+ \sum_{l=l_0}^{l_1-1} V_l \big(S^{\mathcal{I}_{1,l}, \mathcal{I}_{2,l}}U_{l,\rho(\omega_1)}f -A_{l,\omega_2}(U_{l,\rho(\omega_1)}f)\big).
\label{Q6}
\end{eqnarray}
For fixed $\omega_1$ we have 
\begin{eqnarray*}
\lefteqn{\EE_2\left\| S^{\mathcal{I}_{1,l}, \mathcal{I}_{2,l}}U_{l,\rho(\omega_1)}f -A_{l,\omega_2}(U_{l,\rho(\omega_1)}f)\right\|_{L_p^{\mathcal{I}_{1,l}}}^w}
\nonumber\\
&\le&c\EE_2\|A_{l,\omega_2}(0)\|_{L_q^{\mathcal{I}_{1,l}}}^w
\nonumber\\
&&+c\sup_{g\in L_p^{\mathcal{I}_{1,l}\times \mathcal{I}_{2,l}}\setminus\{0\}}\left(\|g\|_{L_p^{\mathcal{I}_{1,l}\times \mathcal{I}_{2,l}}}^{-w}\EE_2\left\|S^{\mathcal{I}_{1,l}, \mathcal{I}_{2,l}}g-A_{l,\omega_2}(g)\right\|_{L_q^{\mathcal{I}_{1,l}}}^w\right)
\|U_{l,\rho(\omega_1)}f\|_{L_p^{\mathcal{I}_{1,l}\times \mathcal{I}_{2,l}}}^w.
\end{eqnarray*}
Taking the expectation with respect to $\omega_1$ and using  \eqref{P8}, this gives
\begin{eqnarray}
\lefteqn{\E\left\|V_l \big(S^{\mathcal{I}_{1,l}, \mathcal{I}_{2,l}}U_{l,\rho(\omega_1)}f -A_{l,\omega_2}(U_{l,\rho(\omega_1)}f)\big)\right\|_{L_q(D_1)}^w}
\notag\\
&\le& c\EE_1\EE_2\left\| S^{\mathcal{I}_{1,l}, \mathcal{I}_{2,l}}U_{l,\rho(\omega_1)}f -A_{l,\omega_2}(U_{l,\rho(\omega_1)}f)\right\|_{L_q(D_1)}^w
\nonumber\\
&\le&c\EE_2\|A_{l,\omega_2}(0)\|_{L_q^{\mathcal{I}_{1,l}}}^w
\notag\\
&&+c\sup_{g\in L_p^{\mathcal{I}_{1,l}\times \mathcal{I}_{2,l}}\setminus\{0\}}\left(\|g\|_{L_p^{\mathcal{I}_{1,l}\times \mathcal{I}_{2,l}}}^{-w}\EE_2\left\|S^{\mathcal{I}_{1,l}, \mathcal{I}_{2,l}}g-A_{l,\omega_2}(g)\right\|_{L_q^{\mathcal{I}_{1,l}}}^w\right)
\EE_1\|U_{l,\rho(\omega_1)}f\|_{L_p^{\mathcal{I}_{1,l}\times \mathcal{I}_{2,l}}}^w.\quad
\label{A7}
\end{eqnarray}
Finally, by \eqref{S3} and \eqref{X1}
\begin{eqnarray}
\EE_1\|U_{l,\rho(\omega_1)}f\|_{L_p^{\mathcal{I}_{1,l}\times \mathcal{I}_{2,l}}}^w
&\le&
\|T_l\|^w\EE_1\|(P_{l+1,\rho(\omega_1)}-P_{l\rho(\omega_1)})f\|_{L_p(D)}^w
\nonumber\\
&\le&
c\,\EE_1\|(P_{l+1,\rho(\omega_1)}-P_{l\rho(\omega_1)})f\|_{L_p(D)}^w.\label{A8}
\end{eqnarray}
Combining \eqref{Q6}--\eqref{A8} gives \eqref{WO1}.

\end{proof}

In the sequel, we will be interested in the following two special cases of the random variable $\rho$ from \eqref{W6}. Firstly, $\rho=\zeta$, where
$\zeta$ is a uniformly distributed on  $[0,1]^d$ random variable taking all values in $[0,1]^d$, and secondly, let $\rho=\zeta_0$, with $\zeta_0$ a random variable with $\zeta_0(\omega_1)\equiv 0$ for all $\omega_1\in\Omega_1$. Thus, $\zeta_0$ is the deterministic case. 
The following result contains the crucial estimates for the discretization technique. 
\begin{proposition}
\label{pro:6}
Let $r,d\in\N$,  $1 \leq p,q \le \infty$, $r/d>1/p-1/q$, $1\le w\le q$, $w<\infty$.  Then there are constants $c_1,c_2>0$ such that for all $l\in \N_0$ 
\begin{eqnarray}
\label{R9}
\sup_{f\in B_{\mathcal{W}_p^r(D)}}(\EE_1\|f-P_{l\zeta}f\|_{L_q(D)}^w)^{1/w}
&\le& c_1 2^{-rl+\left(\frac{1}{p}-\frac{1}{q}\right)_+dl}.
\end{eqnarray}
Moreover, if the embedding condition \eqref{PA2} holds, then
\begin{eqnarray}
\label{X7} 
\sup_{f\in B_{\mathcal{W}_p^r(D)}}\|f-P_{l\zeta_0}f\|_{L_q(D)}
&\le& c_2 2^{-rl+\left(\frac{1}{p}-\frac{1}{q}\right)_+dl}.
\end{eqnarray}
\end{proposition}
\begin{proof}
Relation \eqref{X7} is contained in \cite{Hei04b}, relation (12), 
while \eqref{R9} was shown in \cite{Hei07},  Proposition 1, for $q<\infty$. The estimate \eqref{R9} also holds for $q=\infty$, since in this case $r/d>1/p$, so the embedding condition \eqref{PA2} holds and \cite{Hei04b}, relation (12) gives 
\begin{eqnarray*}
\sup_{f\in B_{\mathcal{W}_p^r(D)}}\|f-P_{l\zeta(\omega_1)}f\|_{L_q(D)}
&\le& c 2^{-rl+\left(\frac{1}{p}-\frac{1}{q}\right)_+dl}
\end{eqnarray*}
for each $\omega_1\in\Omega_1$, with a constant  not depending on $\omega_1$ (the latter following directly from the proof of (12) in \cite{Hei04b}).

\end{proof}
\begin{corollary}
\label{cor:1}
Let $r\in \N$, $1\le p\le \infty$ and $1\le w\le p$, $w<\infty$. Then the following holds for $\rho=\zeta$, and if the embedding condition \eqref{PA2} is satisfied, also for $\rho=\zeta_0$.
\begin{eqnarray*}
\sup_{f\in B_{\mathcal{W}_p^r(D)}}(\EE_1\|U_{l\rho}f\|_{L_p^{\mathcal{I}_{1,l}\times \mathcal{I}_{2,l}}}^w)^{1/w}&\le& c 2^{-rl}.
\end{eqnarray*}
\end{corollary}
\begin{proof}
From \eqref{W8}, \eqref{S3}, \eqref{X1}, and Proposition
\ref{pro:6} we obtain
\begin{eqnarray*}
\sup_{f\in B_{\mathcal{W}_p^r(D)}}(\EE_1\|U_{l\rho}f\|_{L_p^{\mathcal{I}_{1,l}\times \mathcal{I}_{2,l}}}^w)^{1/w}
&\le&
\|T_l\|\sup_{f\in B_{\mathcal{W}_p^r(D)}}(\EE_1\|P_{l+1,\rho}-P_{l\rho}f\|_{L_p(D)}^w)^{1/w}\le c 2^{-rl}.
\end{eqnarray*}
\end{proof}

\begin{lemma}
\label{lem:1}  
Let $r\in \N$, $1 \leq p,q \le \infty$, $1\le w\le p$, $w<\infty$, and assume  \eqref{BB1}. Then the following holds for $\rho=\zeta$ and, if the embedding condition \eqref{PA2} is satisfied, also for $\rho=\zeta_0$.
There is a constant $c>0$ such that for all $l_1\in \N_0$ 
\begin{equation*}
\sup_{f\in B_{\mathcal{W}_p^r(D)}}\left(\EE_1\|Sf-SP_{l_1\rho}f\|_{L_q(D_1)}^w\right)^{1/w}\le c 2^{-rl_1+\left(\frac{1}{p}-\frac{1}{q}\right)_+d_1l_1}
\end{equation*}
\end{lemma}
\begin{proof}
For $p=q$, and therefore also for $q\le p$, this follows directly from \eqref{R9} of Proposition \ref{pro:6}, thus, we can assume $q>p$.  By Corollary \ref{cor:1} and \eqref{P8}, for $l\in\N_0$
\begin{eqnarray}
\label{B3}
\lefteqn{\|V_l S^{\mathcal{I}_{1,l}, \mathcal{I}_{2,l}} U_{l\rho}\|_{L(\mathcal{W}_p^r(D),L_w(\Omega_1,\Prm_1,L_q(D_1))}=\sup_{f\in B_{\mathcal{W}_p^r(D)}}\left(\EE_1\|V_l S^{\mathcal{I}_{1,l}, \mathcal{I}_{2,l}} U_{l\rho}f\|_{L_q(D_1)}^w\right)^{1/w}}
\notag\\
&\le& \|V_l\|\|S^{\mathcal{I}_{1,l}, \mathcal{I}_{2,l}}:L_p^{\mathcal{I}_{1,l}, \mathcal{I}_{2,l}}
\to L_q^{\mathcal{I}_{1,l}}\|\,\sup_{f\in B_{\mathcal{W}_p^r(D)}}\left(\EE_1\|U_{l\rho}f\|_{L_p^{\mathcal{I}_{1,l}\times \mathcal{I}_{2,l}}}^w\right)^{1/w}
\notag\\
&\le&   c2^{-rl+\left(\frac{1}{p}-\frac{1}{q}\right)_+d_1l},
\end{eqnarray}
where we used that 
$$
\big\|S^{\mathcal{I}_{1,l}, \mathcal{I}_{2,l}}:L_p^{\mathcal{I}_{1,l}, \mathcal{I}_{2,l}}
\to L_q^{\mathcal{I}_{1,l}}\big\|\le c2^{\left(\frac{1}{p}-\frac{1}{q}\right)_+d_1l},
$$
see relation (57) 
%   REF-{WN2} 
of Part I.
Since by \eqref{BB1}, $r>\left(\frac{1}{p}-\frac{1}{q}\right)_+d_1$, \eqref{B3} implies for $l_1\in\N_0$
\begin{equation}
\label{AG3}
\sum_{l=l_1}^\infty V_l S^{\mathcal{I}_{1,l}, \mathcal{I}_{2,l}} U_{l\rho}=X
\end{equation}
 with convergence in $L(\mathcal{W}_p^r(D),L_w(\Omega_1,\Prm_1,L_q(D_1))$, and
\begin{equation}
\|X\|_{L(\mathcal{W}_p^r(D),L_w(\Omega_1,\Prm_1,L_q(D_1))}\le c 2^{-rl_1+\left(\frac{1}{p}-\frac{1}{q}\right)_+d_1l_1}.
\label{B4}
\end{equation}
 On the other hand, by \eqref{C1} and \eqref{AG4}, for $m>l_1$ 
$$
SP_{m\rho}-SP_{l_1\rho}=\sum_{l=l_1}^{m-1} V_l S^{\mathcal{I}_{1,l}, \mathcal{I}_{2,l}} U_{l\rho}.
$$
Since $w\le p<q$, Proposition \ref{pro:6} gives
\begin{eqnarray*}
\lefteqn{\sup_{f\in B_{\mathcal{W}_q^r(D)}}\left(\EE_1\bigg\|Sf-SP_{l_1\rho}f-\sum_{l=l_1}^{m} V_l S^{\mathcal{I}_{1,l}, \mathcal{I}_{2,l}} U_{l\rho}f\bigg\|_{L_q(D_1)}^w\right)^{1/w}}
\\
&=& \sup_{f\in B_{\mathcal{W}_q^r(D)}} \left(\EE_1\|Sf-SP_{m\rho}f\|_{L_q(D_1)}^w\right)^{1/w}
\\
&\le& \|S:L_q(D)\to L_q(D_1)\|\,\sup_{f\in B_{\mathcal{W}_q^r(D)}}\left(\EE_1\|f-P_{m\rho}f\|_{L_q(D)}^w\right)^{1/w}
\le c2^{-rm},
\end{eqnarray*}
thus, 
\begin{equation}
\label{AG5}
\sum_{l=l_1}^\infty V_l S^{\mathcal{I}_{1,l}, \mathcal{I}_{2,l}} U_{l\rho}=S-SP_{l_1\rho}
\end{equation}
with convergence in $L(\mathcal{W}_q^r(D),L_w(\Omega_1,\Prm_1,L_q(D_1)))$. Let $J:\mathcal{W}_q^r(D)\to \mathcal{W}_p^r(D)$ denote the identical embedding. Then \eqref{AG3} implies
\begin{equation*}
\sum_{l=l_1}^\infty V_l S^{\mathcal{I}_{1,l}, \mathcal{I}_{2,l}} U_{l\rho}J=XJ
\end{equation*}
with convergence in $L(\mathcal{W}_q^r(D),L_w(\Omega_1,\Prm_1,L_q(D_1))$.
This combined with \eqref{AG5} yields $Xf=Sf-SP_{l_1\rho}f$ for $f\in \mathcal{W}_q^r(D)$ and, since $\mathcal{W}_q^r(D)$ is dense in $\mathcal{W}_p^r(D)$, we have  $X=S-SP_{l_1\rho}$. Now the desired estimate follows from  \eqref{B4}.

\end{proof}
Let us note that if $\rho=\zeta_0$, then $ U_{l,\rho(\omega_1)} f$ does not depend on $\omega_1$ and therefore \eqref{AM5} follows directly from the assumption that the $(A_{l,\omega_2})_{\omega_2\in\Omega_2}$ are randomized algorithms for problem  \eqref{AL8}.
The following is a consequence of Propositions  \ref{pro:5}, \ref{pro:6}, and Lemma \ref{lem:1}.
\begin{corollary}
\label{cor:2}
Let $r\in \N$, $1 \le p,q \le \infty $, $1\le w\le p$, $w<\infty$, let either $\rho=\zeta$ and assume  \eqref{BB1} and \eqref{AM5} or $\rho=\zeta_0$ and assume \eqref{PA2}. Then there is a constant $c>0$ such that for all $l_0,l_1\in\N_0$ with $l_0\le l_1$  
\begin{eqnarray}
	\label{B5}
			\lefteqn{\sup_{f\in B_{\mathcal{W}_p^r(D)}}\left(\E \|Sf-A_{\omega}(f)\|_{L_q(D_1)}^w\right)^{1/w}}
\nonumber\\
		&\le& c2^{-rl_1+\left(\frac{1}{p}-\frac{1}{q}\right)_+d_1l_1}
+c\sum_{l=l_0}^{l_1-1}\left(\EE_2\|A_{l,\omega_2}(0)\|_{L_q^{\mathcal{I}_{1,l}}}^w\right)^{1/w}
\nonumber\\
&&+ c\sum_{l=l_0}^{l_1-1}2^{-rl}\sup_{g\in L_p^{\mathcal{I}_{1,l}\times \mathcal{I}_{2,l}}\setminus\{0\}}\left(\|g\|_{L_p^{\mathcal{I}_{1,l}\times \mathcal{I}_{2,l}}}^{-w}\EE_2\left\|S^{\mathcal{I}_{1,l}, \mathcal{I}_{2,l}}g-A_{l,\omega_2}g\right\|_{L_q^{\mathcal{I}_{1,l}}}^w\right)^{1/w}.
\end{eqnarray}
\end{corollary}

\section{Algorithms for Parametric Integration}
\label{sec:4}
In this and the following sections we assume $d_1,d_2\in\N$. The corresponding results for the borderline cases of integration and approximation can be found in \cite{Hei12}, which also contains references to the vast literature on the complexity of these problems.  Define for $n\in\N$, $n\ge n(0)$
\begin{eqnarray}
&&l_0(n) = d_1\left\lceil \frac{\log n}{d_1(d_1+d_2)} \right\rceil, \quad l_1(n) = \left\lceil \frac{(d_1+d_2)l_0(n)-\sigma_1\log l_0(n) }{d_1}\right\rceil\label{Y1}
\\[.2cm]
&& n_l(n)=\left\lceil  2^{(d_1+d_2)l_0(n) -\delta\min((l-l_0(n)),(l_1(n)-l))}\right\rceil\qquad (l_0(n)\le l<l_1(n)),\label{WO7}
\end{eqnarray}
where the constant $n(0)\in\N$ is chosen in such a way that for each $n\ge n(0)$ the resulting $l_0(n),l_1(n)$ satisfy
\begin{equation}
\label{Z6}
2\le l_0(n)<l_1(n),
\end{equation}
while  the constant $\delta\ge  0$ will be specified later on in each of the considered cases. Note that \eqref{Y1} and \eqref{Z6} imply
\begin{equation}
\label{AK4}
n(0)\ge 2.
\end{equation}
The definition of $l_0(n)$ implies 
\begin{eqnarray}
n^{\frac{1}{d_1+d_2}}&\le& 2^{l_0(n)} \le cn^{\frac{1}{d_1+d_2}},\label{AC7}
\end{eqnarray}
while \eqref{WO7} and \eqref{Z6} give
\begin{equation}
\label{AJ7}
c_1\log(n+1)\le l_0(n)<l_1(n)\le c_2\log(n+1).
\end{equation}
From the definition of $l_1(n)$ we conclude furthermore
\begin{eqnarray}
\label{AJ9}
 2^{\frac{d_1+d_2}{d_1}l_0(n)}l_0(n)^{-\frac{\sigma_1}{d_1}}&\le& 2^{l_1(n)}< 2^{\frac{d_1+d_2}{d_1}l_0(n)+1}l_0(n)^{-\frac{\sigma_1}{d_1}}.
\end{eqnarray}
Combining this with \eqref{AC7} and \eqref{AJ7}, we arrive at 
\begin{eqnarray}
c_3n^{\frac{1}{d_1}}(\log(n+1))^{-\frac{\sigma_1}{d_1}}&\le& 2^{l_1(n)}\le c_4n^{\frac{1}{d_1}}(\log(n+1))^{-\frac{\sigma_1}{d_1}}.\label{AD9}
\end{eqnarray}
Also note that by \eqref{Y1}
\begin{equation}
\label{Z4}
l_1(n)\le \frac{d_1+d_2}{d_1}l_0(n),\quad \text{with equality if  } \sigma_1=0.
\end{equation}

We begin with the non-adaptive setting, where we use Proposition \ref{pro:5}  with 
$$
A_{l,\omega_2}=A_{n_l(n),\omega_2}^{(2)}:L_p^{N_{1,l}, N_{2,l}}\to L_q^{N_{1,l}}\quad(l_0(n)\le l<l_1(n))
$$  
standing for the algorithm from Proposition \ref{pro:4}, $N_{1,l}, N_{2,l}$ as given by \eqref{eq:26}, and we identify $\mathcal{I}_{1,l}$ with $\Z[1,N_{1,l}]$ and $\mathcal{I}_{2,l}$ with $\Z[1,N_{2,l}]$. 
According to \eqref{WM6},  we define
\begin{equation*}
A_{n,\omega}^{(4)}(f)=SP_{l_0(n),\zeta(\omega_1)}f +\sum_{l=l_0(n)}^{l_1(n)-1} V_lA_{n_l(n),\omega_2}^{(2)} ( U_{l,\zeta(\omega_1)} f)\quad(\omega=(\omega_1,\omega_2)).
\end{equation*}
Let us check the measurability condition for $\big(A_{n,\omega}^{(4)}\big)_{\omega\in\Omega}$, that is \eqref{AM5} for $\rho=\zeta$. It follows from the definition of $A_{n_l(n),\omega_2}^{(2)}$ (in Part I) that  the mapping 
\begin{equation}
\label{AM2}
\Omega_2\ni \omega_2 \to A_{n_l(n),\omega_2}^{(2)}\in L\big(L_p^{N_{1,l}, N_{2,l}}, L_q^{N_{1,l}}\big)
\end{equation}
is a random variable taking only a finite number of values and that $\ca\big(A_{n_l(n),\omega_2}^{(2)},g\big)$ neither depends on $\omega_2$ nor on $g$. Consequently, the mappings 
\begin{equation}
\label{AM6}
L_p^{N_{1,l}, N_{2,l}}\times \Omega_2\ni (g,\omega_2) \to \left\{\begin{array}{lll}
   \ca\big(A_{n_l(n),\omega_2}^{(2)},g\big) \\[,2cm]  
A_{n_l(n),\omega_2}^{(2)}(g) 
    \end{array}
\right. 
\end{equation}
are $\mathcal{B}\big(L_p^{N_{1,l}, N_{2,l}}\big)\times\Sigma_2$ measurable. Furthermore, for fixed $f\in \mathcal{W}_p^r(D)$ the mapping $\omega_1\to f(s_{li} +2^{-l}t_{jk}(\zeta(\omega_1)))$ is $\Sigma_1$-measurable. Consequently,  by \eqref{WF3}
\begin{equation}
\label{AM1}
\omega_1\to U_{l\zeta(\omega_1)}f
\end{equation}
is $\Sigma_1$-to-$\mathcal{B}\big(L_p^{N_{1,l}, N_{2,l}}\big)$ measurable.
This implies that \eqref{AM5} holds,
thus $\big(A_{n,\omega}^{(4)}\big)_{\omega\in\Omega}$ is a non-adaptive randomized algorithm.

To state the next result, we set
\begin{eqnarray}
\beta_1&=&\left\{\begin{array}{lll}
	1 & \mbox{if} \quad \frac{r}{d_1}= 1-\frac{1}{\bar{p}}+\left(\frac{1}{p}-\frac{1}{q}\right)_+   \\
	0 & \mbox{otherwise,}    \\
	\end{array}
	\right.
\label{S9} 
\end{eqnarray}
and introduce for $r\in \N_0, d_1,d_2\in\N, 1\le p,q\le \infty$ the following function of $n\in\N$
\begin{eqnarray}
\Phi_1(n)=\left\{ \begin{array}{lll}
n^{\frac{-r+\left(\frac{1}{p}-\frac{1}{q}\right)_+ d_1 -\left(1-\frac{1}{\bar{p}}\right)d_2}{d_1+d_2}}(\log(n+1))^{\sigma_1/2}
& \text{ if }\;  \frac{r}{d_1}> 1-\frac{1}{\bar{p}}+\left(\frac{1}{p}-\frac{1}{q}\right)_+
\\[.3cm] 
n^{-\frac{r}{d_1}+\left(\frac{1}{p}-\frac{1}{q}\right)_+} \,(\log(n+1))^{\sigma_1\left(\frac{r}{d_1}-\left(\frac{1}{p}-\frac{1}{q}\right)_+\right)}
& \text{ if } \; \frac{r}{d_1}\le 1-\frac{1}{\bar{p}}+\left(\frac{1}{p}-\frac{1}{q}\right)_+ ,
\end{array} \right.
\label{T1}
\end{eqnarray}
where we recall that $\sigma_1$ was defined in \eqref{Q8}.
\begin{proposition}
\label{pro:2}
Let $r,d_1,d_2\in\N$, $1 \le p,q \le \infty $, $1\le w\le p$, $w<\infty$, and assume that \eqref{BB1} holds. Then there are constants $c_1,c_2>0$ such that  for all $n\ge n(0)$ and $f\in \mathcal{W}_p^r(D)$
\begin{eqnarray}
	\label{AG9}
			\left(\E \|Sf-A_{n,\omega}^{(4)}(f)\|_{L_q(D_1)}^w\right)^{1/w}
		&\le& c_1\Phi_1(n)(\log (n+1))^{\beta_1}.
\end{eqnarray}
Moreover, 
\begin{equation}
\label{AH0}
\ca(A_{n,\omega}^{(4)},\mathcal{W}_p^r(D))\le c_2n(\log(n+1))^{\beta_1}\quad (\omega\in\Omega).
\end{equation}
\end{proposition}
\begin{proof}
We fix $n\in\N$, $n\ge n(0)$ and write shorthand $l_0$ instead of $l_0(n)$ and respectively $l_1$ and $n_l$.
Taking into account that by \eqref{eq:24} $\kappa'\ge2$, we conclude from \eqref{WO7}, \eqref{AC7}, and \eqref{eq:26} for $l\ge l_0$
\begin{eqnarray*}
n_l&\le& 2^{(d_1+d_2)l_0}<\kappa'2^{(d_1+d_2)l}=N_{1,l}N_{2,l}
\end{eqnarray*}
First we estimate the cardinality of $A_{n,\omega}^{(4)}$. By  \eqref{K1} of Proposition \ref{pro:4} and \eqref{U4} of Proposition \ref{pro:5},
\begin{eqnarray}
\label{WO8}
\ca(A_{n,\omega}^{(4)},B_{ \mathcal{W}_p^r(D)})&\le&  \kappa 2^{(d_1+d_2)l_0}+2\kappa'\kappa''\sum_{l=l_0}^{l_1-1} n_l\nonumber\\
 & \leq &   
\kappa 2^{(d_1+d_2)l_0}+2\kappa'\kappa''\sum_{l=l_0}^{l_1-1} \left(2^{(d_1+d_2)l_0 -\delta\min((l-l_0),(l_1-l))}+1\right)\nonumber\\
& \leq & \left\{\begin{array}{lll}
  c2^{(d_1+d_2)l_0}l_1 & \mbox{if} \quad \delta=0,   \\[.1cm]
 c2^{(d_1+d_2)l_0}  &  \mbox{if} \quad \delta>0.  \\
    \end{array}
\right.
\end{eqnarray}
Furthermore, by \eqref{K2} of Proposition \ref{pro:4},
\begin{eqnarray}
\lefteqn{2^{-rl}\sup_{g\in L_p^{N_{1,l},N_{2,l}}\setminus\{0\}}\left(\|g\|_{L_p^{N_{1,l},N_{2,l}}}^{-w}\EE_2\|S^{N_{1,l}, N_{2,l}}g-A_{n_l,\omega_2}^{(2)}g\|_{L_q^{N_{1,l}}}^w\right)^{1/w}}
\nonumber\\
&\le &
	c2^{-rl}N_{1,l}^{\left(\frac{1}{p}-\frac{1}{q}\right)_+ + 1-\frac{1}{\bar{p}}}	
		n_l^{-\left(1-\frac{1}{\bar{p}}\right)}(\log (N_{1,l}+1))^{\sigma_1/2} 
\le c2^{\gamma(l)+\left(1-\frac{1}{\bar{p}}\right) \delta\min((l-l_0),(l_1-l))}
\label{K3}
\end{eqnarray}
with 
\begin{eqnarray}
\gamma(l)
 &=&-rl +\left(\left(\frac{1}{p}-\frac{1}{q}\right)_+ +1-\frac{1}{\bar{p}}\right)d_1l
-\left(1-\frac{1}{\bar{p}}\right)(d_1+d_2)l_0+\frac{\sigma_1}{2}\log l_0.
\label{X4}
\end{eqnarray}
Note that $\gamma(l)$ considered as a function of the real variable $l$, is linear, with constant (meaning independent of $n$) ascent.  Corollary \ref{cor:2} together with \eqref{AG7} and \eqref{K3} gives
\begin{eqnarray}
\label{X8}
\left(\E \|Sf-A_{n,\omega}^{(4)}(f)\|_{L_q(D_1)}^w\right)^{1/w}
&\le& c2^{-rl_1+\left(\frac{1}{p}-\frac{1}{q}\right)_+d_1l_1}
+c\sum_{l=l_0}^{l_1-1}2^{\gamma (l)+\left(1-\frac{1}{\bar{p}}\right)\delta\min((l-l_0),(l_1-l))}.
\end{eqnarray}
Also observe that 
$$
\gamma (l)+\left(1-\frac{1}{\bar{p}}\right)\delta\min((l-l_0),(l_1-l))
$$
is continuous on $[l_0,l_1]$ and linear on $[l_0,(l_0+l_1)/2]$ and  $[(l_0+l_1)/2,l_1]$, also with constant ascent on each of the intervals. Now we distinguish two cases. 

Case 1: Assume that
\begin{equation}
\label{WQ1}
\frac{r}{d_1}> 1-\frac{1}{\bar{p}}+\left(\frac{1}{p}-\frac{1}{q}\right)_+.
\end{equation}
Then $\gamma(l)$ is a strictly decreasing function of $l$. So we can choose $\delta>0$ in such a way that also
\begin{equation*}
\gamma(l)+\left(1-\frac{1}{\bar{p}}\right)\delta\min((l-l_0),(l_1-l)) \quad\text{is a strictly decreasing function of } l.
\end{equation*}
Consequently,
\begin{eqnarray}
\label{WQ4}
\sum_{l=l_0}^{l_1-1}2^{\gamma (l)+\left(1-\frac{1}{\bar{p}}\right)\delta\min((l-l_0),(l_1-l))}\le c 2^{\gamma(l_0)}.
\end{eqnarray}
We have 
\begin{eqnarray}
\label{X9}
\gamma(l_0)
 &=&\left(-r +\left(\frac{1}{p}-\frac{1}{q}\right)_+d_1 -\left(1-\frac{1}{\bar{p}}\right)d_2\right) l_0 +\frac{\sigma_1}{2}\log l_0.
\end{eqnarray}
Note that by \eqref{Y1}, \eqref{WQ1}, and \eqref{X9}
\begin{eqnarray*}
\lefteqn{\left(-r+\left(\frac{1}{p}-\frac{1}{q}\right)_+d_1\right)l_1\le-\left(r-\left(\frac{1}{p}-\frac{1}{q}\right)_+d_1\right)\frac{(d_1+d_2)l_0 -\sigma_1\log  l_0}{d_1}+c}
\nonumber\\
&\le& -\left(r-\left(\frac{1}{p}-\frac{1}{q}\right)_+d_1\right)l_0-\left(r-\left(\frac{1}{p}-\frac{1}{q}\right)_+d_1\right)\frac{d_2l_0 -\sigma_1\log  l_0}{d_1}+c
\nonumber\\
&\le& -\left(r-\left(\frac{1}{p}-\frac{1}{q}\right)_+d_1\right)l_0-\left(1-\frac{1}{\bar{p}}\right)\left(d_2 l_0 -\sigma_1\log  l_0\right) +c\le \gamma(l_0)+c,
\end{eqnarray*}
which together with \eqref{X8} and \eqref{WQ4} yields 
\begin{eqnarray*}
\lefteqn{\sup_{f\in B_{ \mathcal{W}_p^r(D)}} \left(\E \|Sf-A^{(4)}_{n,\omega}f\|_{L_p(D_1)}^w\right)^{1/w}}
\notag\\
&\le&c  2^{\gamma(l_0)}\le c  2^{\left(-r+\left(\frac{1}{p}-\frac{1}{q}\right)_+d_1-\left(1-\frac{1}{\bar{p}}\right)d_2\right)l_0+\frac{\sigma_1}{2}\log l_0}
\le cn^{-\frac{r-\left(\frac{1}{p}-\frac{1}{q}\right)_+d_1+\left(1-\frac{1}{\bar{p}}\right)d_2}{d_1+d_2}}(\log(n+1))^{\frac{\sigma_1}{2}},
\end{eqnarray*}
the last relation being a consquence of \eqref{AC7} and \eqref{AJ7}.
Moreover, (\ref{WO8}) and and \eqref{AC7} imply $\ca(A^{(4)}_{n,\omega},B_{ \mathcal{W}_p^r(D)})\le cn$. This proves \eqref{AG9} and \eqref{AH0}  in the first case.

Case 2: Now let
\begin{equation}
\label{WQ2}
 \frac{r}{d_1}\le 1-\frac{1}{\bar{p}}+\left(\frac{1}{p}-\frac{1}{q}\right)_+.
\end{equation}
If equality holds in \eqref{WQ2}, or equivalently, $\gamma(l_0)= \gamma(l_1)$, we set $\delta=0$, otherwise  by \eqref{X4}, $\gamma(l) $ is a strictly increasing function of  $l$.
 Here we choose   $\delta>0$ in such a way that also
\begin{equation*}
\gamma(l)+\left(1-\frac{1}{\bar{p}}\right)\delta\min((l-l_0),(l_1-l))  \quad\text{is a strictly increasing function of } l.
\end{equation*}
 This implies
\begin{eqnarray}
\label{Y0}
\sum_{l=l_0}^{l_1-1}2^{\gamma (l)+\left(1-\frac{1}{\bar{p}}\right) \delta\min((l-l_0),(l_1-l))}
\le cl_1^{\beta_1} 2^{\gamma (l_1)}. 
\end{eqnarray}
Relation \eqref{X4} together with \eqref{Y1} 
gives
\begin{eqnarray}
\label{Y2}
\gamma(l_1)&=&-rl_1+\left(\frac{1}{p}-\frac{1}{q}\right)_+d_1l_1
\nonumber\\
&&+\left(1-\frac{1}{\bar{p}}\right)d_1l_1
-\left(1-\frac{1}{\bar{p}}\right)(d_1+d_2)l_0+\frac{\sigma_1}{2}\log l_0
\nonumber\\
&\le&-rl_1+\left(\frac{1}{p}-\frac{1}{q}\right)_+d_1l_1+c,
\end{eqnarray}
where we used that $\sigma_1=1$ implies $\bar{p}=2$.
From \eqref{X8}, \eqref{Y0}, and \eqref{Y2} we obtain
\begin{eqnarray}
\lefteqn{\sup_{f\in B_{ \mathcal{W}_p^r(D)}} \left(\E \|Sf-A^{(4)}_{n,\omega}f\|_{L_p(D_1)}^w\right)^{1/w} }\nonumber
\\ 
& \le &
c\,2^{\left(-r+\left(\frac{1}{p}-\frac{1}{q}\right)_+d_1\right)l_1} +cl_1^{\beta_1} 2^{\gamma (l_1)}.
 \le cl_1^{\beta_1}2^{\left(-r+\left(\frac{1}{p}-\frac{1}{q}\right)_+d_1\right)l_1}
\nonumber\\
&\le& cn^{-\frac{r}{d_1}+\left(\frac{1}{p}-\frac{1}{q}\right)_+} (\log(n+1))^{\sigma_1\left(\frac{r}{d_1}-\left(\frac{1}{p}-\frac{1}{q}\right)_+\right)+\beta_1},
\label{AC3}
\end{eqnarray}
where we also used \eqref{AJ7} and \eqref{AD9}.
This proves \eqref{AG9} in the second case.

If $\frac{r}{d_1}< 1-\frac{1}{\bar{p}}+\left(\frac{1}{p}-\frac{1}{q}\right)_+$, we have $\beta_1=0$, $\delta>0$, so  by \eqref{WO8} and \eqref{AC7}, $\ca(A_{n,\omega}^{(4)},B_{ \mathcal{W}_p^r(D)})\le cn$, while 
if $\frac{r}{d_1}= 1-\frac{1}{\bar{p}}+\left(\frac{1}{p}-\frac{1}{q}\right)_+$, we have 
 $\beta_1=1$, $\delta=0$, so (\ref{WO8}) together with \eqref{AJ7} and \eqref{AD9}  gives $\ca(A_{n,\omega}^{(4)},B_{ \mathcal{W}_p^r(D)})\le cn\log(n+1)$, which proves
\eqref{AH0} also in the second case.

\end{proof}

Now we turn to the adaptive case and $2<p<q$. For $n\in\N$, $n\ge n(0)$ let $l_0(n)$, $l_1(n)$, and $n_l(n)$ be given by \eqref{Y1} and \eqref{WO7}.
Here we define for $l_0(n)\le l<l_1(n)$
$$
m_l=\lceil c(1)\log(N_{1,l}+N_{2,l})\rceil,
$$
where $c(1)$ stands for the constant $c_1$ from Proposition \ref{pro:7} and $N_{1,l}, N_{2,l}$ from \eqref{eq:26}. It follows that
\begin{equation*}
m_l\le c l.
\end{equation*}
We use  \eqref{WM6} with  
$$
A_{l,\omega_2}=A_{n_l(n),m_l,\omega_2}^{(3)}:L_p^{N_{1,l}, N_{2,l}}\to L_q^{N_{1,l}}
$$
the algorithm from Proposition \ref{pro:7} and define an adaptive algorithm
\begin{equation*}
A_{n,\omega}^{(5)}(f)=SP_{l_0(n),\zeta(\omega_1)}f +\sum_{l=l_0(n)}^{l_1(n)-1} V_lA_{n_l(n),m_l,\omega_2}^{(3)} ( U_{l,\zeta(\omega_1)} f)\quad(\omega=(\omega_1,\omega_2)).
\end{equation*}
Similarly to the non-adaptive case \eqref{AM2}--\eqref{AM1} let us verify the measurability condition for $\big(A_{n,\omega}^{(5)}\big)_{\omega\in\Omega}$. Here we note that the definition of $A_{n_l(n),m_l,\omega_2}^{(3)}$ in Part I implies that 
\begin{equation}
\label{AM2a}
\Omega_2\ni \omega_2 \to A_{n_l(n),m_l,\omega_2}^{(3)}\in \mathcal{F}\big(L_p^{N_{1,l}, N_{2,l}}, L_q^{N_{1,l}}\big)
\end{equation}
is a random variable taking only a finite number of values and that for each $\omega_2$ the  mappings 
\begin{equation}
\label{AM3b}
L_p^{N_{1,l}, N_{2,l}}\ni g \to \left\{\begin{array}{lll}
 \ca(A_{n_l(n),m_l,\omega_2}^{(3)},g)    \\[.2cm]
  A_{n_l(n),m_l,\omega_2}^{(3)}(g)    
    \end{array}
\right. 
\end{equation}
are $\mathcal{B}\big(L_p^{N_{1,l}, N_{2,l}}\big)$-measurable.   Consequently, the mappings 
\begin{equation}
\label{AM3a}
L_p^{N_{1,l}, N_{2,l}}\times \Omega_2\ni (g,\omega_2) \to \left\{\begin{array}{lll}
 \ca(A_{n_l(n),m_l,\omega_2}^{(3)},g)    \\[.2cm]
  A_{n_l(n),m_l,\omega_2}^{(3)}(g)    
    \end{array}
\right. 
\end{equation}
are $\mathcal{B}\big(L_p^{N_{1,l}, N_{2,l}}\big)\times\Sigma_2$-measurable. Now  fix $f\in \mathcal{W}_p^r(D)$. As already mentioned, see \eqref{AM1},
the mapping $\omega_1\to U_{l\zeta(\omega_1)}f$ 
is $\Sigma_1$-measurable,
implying that \eqref{AM5} is satisfied, hence
  $\big(A_{n,\omega}^{(5)}\big)_{\omega\in\Omega}$ is an adaptive randomized algorithm.

For  $r\in \N_0, d_1,d_2\in\N, 2< p<q\le \infty$ we define the function $\Phi_2(n)$ of $n\in\N$ as follows. If
\begin{equation}
\label{B7}
\bigg(\frac{1}{2}-\frac{1}{p}\bigg)d_2 > \bigg(\frac{1}{p}-\frac{1}{q}\bigg)d_1,
\end{equation}
then we set
\begin{eqnarray}
\Phi_2(n)&=&
 \left\{ \begin{array}{lll}
n^{\frac{-r-\frac12 d_2}{d_1+d_2}}
& \text{ if }\;  \frac{r}{d_1}>  \left(\frac{1}{p}-\frac{1}{q}\right)\left(\frac{d_1}{d_2}+1\right)+\frac{1}{2}
\\[.3cm] 
n^{-\frac{r}{d_1}+\frac{1}{p}-\frac{1}{q}} 
& \text{ if } \;  \frac{r}{d_1}\le \left(\frac{1}{p}-\frac{1}{q}\right)\left(\frac{d_1}{d_2}+1\right)+\frac{1}{2},
\end{array} \right.
\label{C7}
\end{eqnarray}
while if
\begin{equation}
\label{B8}
\bigg(\frac{1}{2}-\frac{1}{p}\bigg)d_2\le \bigg(\frac{1}{p}-\frac{1}{q}\bigg)d_1,
\end{equation}
we define
\begin{eqnarray}
\Phi_2(n)&=& \left\{ \begin{array}{lll}
n^{\frac{-r+\big(\frac{1}{p}-\frac{1}{q}\big)d_1-\big(1-\frac{1}{p}\big)d_2 }{d_1+d_2}}
& \text{ if }\;  \frac{r}{d_1}> 1-\frac{1}{q}
\\[.3cm] 
n^{-\frac{r}{d_1}+\frac{1}{p}-\frac{1}{q}} 
& \text{ if } \;  \frac{r}{d_1}\le 1-\frac{1}{q}.
\end{array} \right.
\label{C0}
\end{eqnarray}

The following observations, which are easily checked,  complement the definition of $\Phi_2$ and will be of help in the sequel: 
\begin{eqnarray}
&&\frac{-r-\frac12 d_2}{d_1+d_2}\ge \frac{-r+\big(\frac{1}{p}-\frac{1}{q}\big)d_1-\big(1-\frac{1}{p}\big)d_2 }{d_1+d_2}\quad
\text{iff}\quad 1-\frac{1}{q}  \ge \left(\frac{1}{p}-\frac{1}{q}\right)\left(\frac{d_1}{d_2}+1\right)+\frac{1}{2}
\notag\\
&&\text{iff}\quad\left(\frac{1}{2}-\frac{1}{p}\right)d_2\ge\left(\frac{1}{p}-\frac{1}{q}\right)d_1
\label{C4A}
\\[.3cm]
&&\frac{-r-\frac12 d_2}{d_1+d_2}\ge-\frac{r}{d_1}+\frac{1}{p}-\frac{1}{q}
\quad 
\text{iff}\quad \frac{r}{d_1}\ge \left(\frac{1}{p}-\frac{1}{q}\right)\left(\frac{d_1}{d_2}+1\right)+\frac{1}{2}\quad
\label{C5A}\\[.3cm]
&&\frac{-r+\big(\frac{1}{p}-\frac{1}{q}\big)d_1-\big(1-\frac{1}{p}\big)d_2 }{d_1+d_2}\ge-\frac{r}{d_1}+\frac{1}{p}-\frac{1}{q}\quad
\text{iff}\quad \frac{r}{d_1}\ge 1-\frac{1}{q},
\label{C6A}
\end{eqnarray}
where, moreover, in each of \eqref{C4A}, \eqref{C5A}, and \eqref{C6A} equality holds in one relation iff it holds all relations.
This implies, in particular, that the passage from one rate to another in \eqref{C7} and \eqref{C0} is ''continuous'', or in other words, we may replace ''$>$'' by ''$\ge$'' in \eqref{B7}, \eqref{C7}, and \eqref{C0} without producing a contradiction.

Finally we put
\begin{eqnarray}
\beta_2&=&\left\{\begin{array}{lll}
	1 & \mbox{if} \quad \Big(\frac{1}{2}-\frac{1}{p}\Big)d_2\le \Big(\frac{1}{p}-\frac{1}{q}\Big)d_1\quad \mbox{and}\quad\frac{r}{d_1}= 1-\frac{1}{q}   \\
	0 & \mbox{otherwise}    \\
	\end{array}
	\right.
\label{Y5} 
\end{eqnarray}
and recall that $n(0)$ was defined before \eqref{Z6}.
\begin{proposition}
\label{pro:3}
Let $r,d_1,d_2\in\N$, $2<p< q \le \infty$, $1\le w\le p$, $w<\infty$, and assume that \eqref{BB1} holds. Then there are constants $c_1,c_2>0$ such that  for all $n\ge n(0)$ and $f\in \mathcal{W}_p^r(D)$
\begin{eqnarray}
	\label{AH1}
			\left(\E \|Sf-A_{n,\omega}^{(5)}(f)\|_{L_q(D_1)}^w\right)^{1/w}
		&\le& c_1\Phi_2\left(n\right)(\log(n+1))^{\beta_2}.
\end{eqnarray}
Moreover, 
\begin{equation}
\label{AH2}
\ca(A_{n,\omega}^{(5)},\mathcal{W}_p^r(D))\le c_2n(\log(n+1))^{\beta_2+1}.
\end{equation}
\end{proposition}
\begin{proof}
By \eqref{WM1} of Proposition \ref{pro:7} and \eqref{U4} of Proposition \ref{pro:5}, the cardinality of $A^{(5)}_{n,\omega}$ satisfies. 
\begin{eqnarray}
\label{WO9}
\ca(A_{n,\omega}^{(5)},B_{ \mathcal{W}_p^r(D)})&\le&  \kappa 2^{(d_1+d_2)l_0}+6\kappa'\kappa''\sum_{l=l_0}^{l_1-1} m_ln_l\nonumber\\
 & \leq &   
c 2^{(d_1+d_2)l_0}+cl_1\sum_{l=l_0}^{l_1-1} \left(2^{(d_1+d_2)l_0 -\delta\min((l-l_0),(l_1-l))}+1\right) \nonumber\\
& \leq & \left\{\begin{array}{lll}
  c2^{(d_1+d_2)l_0}l_1^2 & \mbox{if} \quad \delta=0,   \\[.1cm]
 c2^{(d_1+d_2)l_0} l_1 &  \mbox{if} \quad \delta>0.  
    \end{array}
\right.
\end{eqnarray}
Moreover,  by  \eqref{mc-eq:2} of Proposition \ref{pro:7}, 
\begin{eqnarray}
	\label{K4}
	\lefteqn{2^{-rl}\sup_{g\in L_p^{N_{1,l},N_{2,l}}\setminus\{0\}}\left(\|g\|_{L_p^{N_{1,l},N_{2,l}}}^{-w}\EE_2\|S^{N_{1,l}, N_{2,l}}g-A_{n_l,m_l\omega_2}^{(3)}g\|_{L_q^{N_{1,l}}}^w\right)^{1/w}}
\nonumber\\
&\le&c2^{-rl}\Bigg(N_{1,l}^{1/p-1/q}\left\lceil\frac{n_l}{N_{1,l}}\right\rceil^{-\left(1-1/p\right)}+\left\lceil\frac{n_l}{N_{1,l}}\right\rceil^{-1/2}\Bigg)\|f\|_{L_p^{N_{1,l},N_{2,l}}}
\nonumber\\
&\le &
	c2^{-rl}\left(N_{1,l}^{1-1/q}	
		n_l^{-1+1/p}+N_{1,l}^{1/2}n_l^{-1/2}\right) 
\nonumber\\
&\le& c2^{\gamma_1(l)+\left(1-\frac{1}{p}\right)\delta\min((l-l_0),(l_1-l))}
 + c2^{\gamma_2(l)+\frac{1}{2}\delta\min((l-l_0),(l_1-l))},
	\end{eqnarray}
where we defined
\begin{eqnarray*}
\gamma_1(l)&=&\left(-r+\left(1-\frac{1}{q}\right)d_1\right)l-\left(1-\frac{1}{p}\right)(d_1+d_2)l_0
\\
\gamma_2(l)&=&\left(-r+\frac{1}{2}d_1\right)l-\frac{1}{2}(d_1+d_2)l_0.
\end{eqnarray*}
Note that we have $\sigma_1=0$ and due to \eqref{Z4},  
$l_1 = \frac{(d_1+d_2)l_0 }{d_1}$. 
Therefore the values at the endpoints of the interval $[l_0,l_1]$ are
\begin{eqnarray}
\gamma_1(l_0)&=& \left(-r+\left(\frac{1}{p}-\frac{1}{q}\right)d_1-\left(1-\frac{1}{p}\right)d_2\right)l_0
\label{AH5}\\
\gamma_1(l_1)&=&
\left(-r\frac{d_1+d_2}{d_1}+\left(\frac{1}{p}-\frac{1}{q}\right)(d_1+d_2)\right)l_0
\label{AH6}\\
\gamma_2(l_0)&=&\left(-r-\frac{d_2}{2}\right)l_0
\label{AH7}\\
\gamma_2(l_1)&=&-r\frac{d_1+d_2}{d_1}l_0.
\label{AH8}
\end{eqnarray}
Observe that $\gamma_1(l_0),\gamma_1(l_1),\gamma_2(l_0),\gamma_2(l_1)$ are constant multiples of $l_0$ and, since by \eqref{Z6} $l_0\ne 0$, we have
\begin{equation*}
\gamma_1(l_1)>\gamma_2(l_1). 
\end{equation*}
It follows from \eqref{C4A}--\eqref{C6A} that
\begin{eqnarray*}
\gamma_2(l_0)&\ge& \gamma_1(l_0)\quad
\text{iff}\quad\left(\frac{1}{2}-\frac{1}{p}\right)d_2\ge \left(\frac{1}{p}-\frac{1}{q}\right)d_1
\\
\gamma_2(l_0)&\ge&\gamma_1(l_1)
\quad 
\text{iff}\quad \frac{r}{d_1}\ge \left(\frac{1}{p}-\frac{1}{q}\right)\left(\frac{d_1}{d_2}+1\right)+\frac{1}{2}
\\
\gamma_1(l_0)&\ge&\gamma_1(l_1)\quad
\text{iff}\quad \frac{r}{d_1}\ge1-\frac{1}{q},
\end{eqnarray*}
where, moreover, in each line equality holds in the left-hand relation iff it holds in the right-hand relation.
We obtain from  Corollary \ref{cor:2} together with \eqref{AG7}, \eqref{K4},  and \eqref{AH6}
\begin{eqnarray}
\label{A1}
\lefteqn{\sup_{f\in B_{ \mathcal{W}_p^r(D)}} \left(\E \|Sf-A_{n,\omega}^{(5)}(f)\|_{L_p(D_1)}^w\right)^{1/w}} 
\nonumber\\ 
& \leq &
c\,2^{\gamma_1(l_1)} +c  \sum_{l=l_0}^{l_1-1}2^{\gamma_1(l)+\left(1-\frac{1}{p}\right)\delta\min((l-l_0),(l_1-l))}
 +c  \sum_{l=l_0}^{l_1-1}2^{\gamma_2(l)+\frac{1}{2}\delta\min((l-l_0),(l_1-l))}.\quad
\end{eqnarray}

Case 1: Suppose that
\begin{equation}
\label{AE7}
\left(\frac{1}{2}-\frac{1}{p}\right)d_2>\left(\frac{1}{p}-\frac{1}{q}\right)d_1,
\end{equation}
hence $\gamma_2(l_0)>\gamma_1(l_0)$. 

Case 1.1: If, moreover,
\begin{equation}
\frac{r}{d_1}> \left(\frac{1}{p}-\frac{1}{q}\right)\left(\frac{d_1}{d_2}+1\right)+\frac{1}{2},
\label{AL1}
\end{equation}
then $\gamma_2(l_0)>\gamma_1(l_1)$, thus summarizing, we have
$$
\gamma_2(l_0)>\gamma_1(l_0),\quad \gamma_2(l_0)>\gamma_1(l_1)>\gamma_2(l_1).
$$  
It follows that $\gamma_2(l)$ is a strictly decreasing function of $l$ and 
$$
\gamma_2(l_0)-\max(\gamma_1(l_0),\gamma_1(l_1))>cl_0
$$  
for some constant $c>0$, so we can choose $\delta>0$ in such a way that also
\begin{equation*}
\gamma_2(l)+\frac{1}{2}\delta\min((l-l_0),(l_1-l)) \quad\text{is a strictly decreasing function of } l
\end{equation*}
and
\begin{equation*}
\gamma_1(l)+\left(1-\frac{1}{p}\right)\delta\min((l-l_0),(l_1-l))\le \gamma_2(l_0)-\delta l_0\quad  (l_0\le l\le l_1).
\end{equation*}
We obtain from \eqref{A1} 
\begin{eqnarray*}
\sup_{f\in B_{ \mathcal{W}_p^r(D)}} \left(\E \|Sf-A_{n,\omega}^{(5)}(f)\|_{L_p(D_1)}^w\right)^{1/w} 
& \leq &
c\,  2^{\gamma_2(l_0)}
\le cn^{\frac{-r-\frac{d_2}{2}}{d_1+d_2}}. 
\end{eqnarray*}

Case 1.2: On the other hand, if \eqref{AE7} holds and  
$$
\frac{r}{d_1}\le  \left(\frac{1}{p}-\frac{1}{q}\right)\left(\frac{d_1}{d_2}+1\right)+\frac{1}{2},
$$
then $\gamma_2(l_0)\le \gamma_1(l_1)$, hence 
$$
\gamma_1(l_1)\ge \gamma_2(l_0)>\gamma_1(l_0),\quad \gamma_1(l_1)>\gamma_2(l_1).
$$ 
If $\gamma_1(l_1)= \gamma_2(l_0)$, then  $\gamma_1(l_1)>\gamma_1(l_0)$ and $\gamma_2(l_0)>\gamma_2(l_1)$, and we choose $\delta>0$  so that  
\begin{eqnarray}
&&\gamma_1(l)+\left(1-\frac{1}{p}\right)\delta\min((l-l_0),(l_1-l)) \quad\text{is a strictly increasing function of } l
\label{Y7}\\
&&\gamma_2(l)+\frac{1}{2}\delta\min((l-l_0),(l_1-l)) \quad\text{is a strictly decreasing function of } l.
\label{Y8}
\end{eqnarray}
If $\gamma_1(l_1)> \gamma_2(l_0)$, then we define $\delta>0$ in such a way that \eqref{Y7} holds and
\begin{eqnarray*}
\gamma_2(l)+\frac{1}{2}\delta\min((l-l_0),(l_1-l)) \le \gamma_1(l_1)-\delta l_0\quad(l_0\le l\le l_1).
\end{eqnarray*}
In each case we obtain from \eqref{A1}
\begin{eqnarray*}
\sup_{f\in B_{ \mathcal{W}_p^r(D)}} \left(\E \|Sf-A_{n,\omega}^{(5)}(f)\|_{L_p(D_1)}^w\right)^{1/w} 
& \leq &
c\,2^{\gamma_1(l_1)}\le cn^{-\frac{r}{d_1}+\frac{1}{p}-\frac{1}{q}}. 
\end{eqnarray*}

Case 2: Now assume
\begin{equation}
\label{AF0}
\left(\frac{1}{2}-\frac{1}{p}\right)d_2\le\left(\frac{1}{p}-\frac{1}{q}\right)d_1,
\end{equation}
hence $\gamma_1(l_0)\ge \gamma_2(l_0)$. 

Case 2.1: If, furthermore, 
$$
\frac{r}{d_1}> 1-\frac{1}{q},
$$
then $\gamma_1(l_0)>\gamma_1(l_1)$,
so altogether
$$
\gamma_1(l_0)\ge \gamma_2(l_0), \quad \gamma_1(l_0)>\gamma_1(l_1)>\gamma_2(l_1).
$$
If $ \gamma_1(l_0)= \gamma_2(l_0)$, then $\gamma_2(l_0)>\gamma_2(l_1)$, meaning that both $ \gamma_1(l)$ and $\gamma_2(l)$ are strictly decreasing, and we choose $\delta>0$ so that 
\begin{eqnarray}
&&\gamma_1(l)+\left(1-\frac{1}{p}\right)\delta\min((l-l_0),(l_1-l)) \quad\text{is a strictly decreasing function of } l
\label{Z0}\\
&&\gamma_2(l)+\frac{1}{2}\delta\min((l-l_0),(l_1-l)) \quad\text{is a strictly decreasing function of } l.
\nonumber
\end{eqnarray}
If $ \gamma_1(l_0)> \gamma_2(l_0)$, we choose $\delta>0$ in such a way that \eqref{Z0} holds and
\begin{eqnarray*}
\gamma_2(l)+\frac{1}{2}\delta\min((l-l_0),(l_1-l)) \le \gamma_1(l_0)-\delta l_0\quad  (l_0\le l\le l_1).
\end{eqnarray*}
We conclude
\begin{eqnarray*}
\sup_{f\in B_{ \mathcal{W}_p^r(D)}} \left(\E \|Sf-A_{n,\omega}^{(5)}(f)\|_{L_p(D_1)}^w\right)^{1/w} 
& \leq &
c\, 2^{\gamma_1(l_0)}\le cn^{\frac{-r+\big(\frac{1}{p}-\frac{1}{q}\big)d_1-\big(1-\frac{1}{p}\big)d_2}{d_1+d_2}}.
\end{eqnarray*}

Case 2.2: On the other hand, if \eqref{AF0} holds and
\begin{equation}
\label{AC6}
\frac{r}{d_1}\le 1-\frac{1}{q},
\end{equation}
then $\gamma_1(l_0)\le \gamma_1(l_1)$,
thus 
$$
\gamma_1(l_1)\ge\gamma_1(l_0)\ge \gamma_2(l_0), \quad \gamma_1(l_1)>\gamma_2(l_1).
$$
If $\gamma_1(l_1)=\gamma_1(l_0)$, then equality holds in \eqref{AC6} and we have $\beta_2=1$, so we set $\delta=0$. 
If $\gamma_1(l_1)>\gamma_1(l_0)$, we choose $\delta>0$ so that 
\begin{eqnarray*}
&&\gamma_1(l)+\left(1-\frac{1}{p}\right)\delta\min((l-l_0),(l_1-l)) \quad\text{is a strictly increasing function of } l,
\\
&&\gamma_2(l)+\frac{1}{2}\delta\min((l-l_0),(l_1-l)) \le \gamma_1(l_1)-\delta l_0\quad  (l_0\le l\le l_1).
\end{eqnarray*}
In both cases \eqref{A1} gives
\begin{eqnarray*}
\sup_{f\in B_{ \mathcal{W}_p^r(D)}} \left(\E \|Sf-A_{n,\omega}^{(5)}(f)\|_{L_p(D_1)}^w\right)^{1/w} 
& \leq &
cl_1^{\beta_2}2^{\gamma_1(l_1)} 
\le cn^{-\frac{r}{d_1}+\frac{1}{p}-\frac{1}{q}} (\log(n+1))^{\beta_2}.
\end{eqnarray*}
Finally note that $\delta=0$ only in the case 2.2 with $\gamma_1(l_1)=\gamma_1(l_0)$, or equivalently, if $\beta_2=1$.  Therefore by \eqref{WO9}
\begin{eqnarray*}
\ca(A_{n,\omega}^{(5)},B_{ \mathcal{W}_p^r(D)})&\le&  
  cn(\log(n+1))^{\beta_2+1},
\end{eqnarray*}
which is \eqref{AH2}.

\end{proof}

\section{Lower Bounds and Complexity}
\label{sec:5}
The following theorem extends Wiegand's Theorem 5.1 of  \cite{Wie06} to  the case of $p\ne q$. Furthermore, Wiegand considered the randomized setting only for the case that  $W_p^r(D)$ is embedded into $C(D)$. We provide an analysis of the  randomized setting for the case  of non-embedding, as well. Recall that the functions $\Phi_1,\Phi_2$ where defined in \eqref{T1} and \eqref{B7}--\eqref{C0}, and the parameters $\beta_1,\beta_2$ in \eqref{S9} and \eqref{Y5}, respectively.
Furthermore, put
\begin{eqnarray}
\sigma_2&=&\left\{\begin{array}{lll}
	1 & \text{if}\;  \left(\Big(\frac{1}{2}-\frac{1}{p}\Big)d_2\ge \Big(\frac{1}{p}-\frac{1}{q}\Big)d_1\right)\wedge \left(\frac{r}{d_1}\ge  \left(\frac{1}{p}-\frac{1}{q}\right)\left(\frac{d_1}{d_2}+1\right)+\frac{1}{2}\right)\wedge \left(q=\infty\right) \\
	0 & \mbox{otherwise.}    \\
	\end{array}
	\right.\quad\label{AE1}
\end{eqnarray}
\begin{theorem}
\label{theo:1}
Let $d_1,d_2\in\N$, $r\in \N_0$, $1 \leq p,q \le \infty$ and assume that \eqref{B1} is satisfied.  Then there are constants $c_{1-8}>0$ such that the following hold for all $n\in \N$. If   $(p\le 2)\vee (p\ge q)$, then 
\begin{eqnarray}
\label{mc-eq:2c} 
c_1 \Phi_1(n)&\le& e_n^\ran(S,B_{\mathcal{W}_p^r(D)},L_q(D_1))
\notag\\
&\le& e_n^\ranno(S,B_{\mathcal{W}_p^r(D)},L_q(D_1))\le c_2 \Phi_1(n)(\log (n+1))^{\beta_1\left(2-\frac{1}{\bar{p}}\right)}.\quad
\end{eqnarray}
On the other hand, if $2<p<q$, then in the non-adaptive setting
\begin{eqnarray}
\label{B6} 
c_3 \Phi_1(n)&\le& e_n^\ranno(S,B_{\mathcal{W}_p^r(D)},L_q(D_1))\le c_4 \Phi_1(n)(\log (n+1))^{\beta_1\left(2-\frac{1}{p}\right)},\quad
\end{eqnarray}
while in the adaptive setting
\begin{eqnarray}
c_5\Phi_2(n)(\log(n+1))^{\sigma_2/2}&\le&  e_n^{\rm ran}(S,B_{\mathcal{W}_p^r(D)},L_q(D_1))
\notag\\
&\le& c_6\Phi_2\left(\frac{n}{\log(n+1)}\right)(\log(n+1))^{\beta_2\left(2-\frac{1}{p}\right)}.
\label{B9}
\end{eqnarray}
Finally, if the embedding condition \eqref{PA2} is satisfied, then in the deterministic setting
\begin{equation}
 c_7n^{\frac{-r+d_1\left(\frac{1}{p}-\frac{1}{q}\right)_+}{d_1+d_2}}\le e_n^\de(S,B_ {\mathcal{W}_p^r(D)},L_q(D_1))\le c_8 n^{\frac{-r+d_1\left(\frac{1}{p}-\frac{1}{q}\right)_+}{d_1+d_2}}.
\label{AE2}
\end{equation}
\end{theorem}

\begin{proof} 
For $n<n(0)$ the upper bounds are trivial, while 
If $\frac{r}{d_1}= \left(\frac{1}{p}-\frac{1}{q}\right)_+$, we have by \eqref{T1} $\Phi_1(n)=1$ and for $2<p<q$ by \eqref{C7} and \eqref{C0} also $\Phi_2(n)=1$ for all $n\in\N$, hence all upper bounds are trivial, so we assume for the proof of the upper bounds $\frac{r}{d_1}> \left(\frac{1}{p}-\frac{1}{q}\right)_+$ and in particular $r\ge 1$.  
 
The upper bound for the deterministic setting follows directly from Corollary \ref{cor:2} with $\rho=\zeta_0$ and $l_1=l_0=\left\lceil \frac{\log n}{d_1+d_2} \right\rceil$, thus, the algorithm in \eqref{WM6} consists only of  the deterministic $SP_{l_0,\zeta_0}$.

To show the upper bounds in the randomized non-adaptive setting, we conclude from  \eqref{AG9} and \eqref{AH0} of Proposition \ref{pro:2} that
\begin{equation}
\label{AH4}
e_{\left\lceil cn(\log(n+1))^{\beta_1}\right\rceil}^\ranno(S,B_{\mathcal{W}_p^r(D)},L_q(D_1))\le   c\Phi_1(n) (\log(n+1))^{\beta_1}.
\end{equation}
This directly yields \eqref{mc-eq:2c} and \eqref{B6} for the case $\beta_1=0$. 
If $\beta_1=1$, we have by \eqref{S9},
\begin{equation}
\label{AH9}
\frac{r}{d_1}=1-\frac{1}{\bar{p}}+\left(\frac{1}{p}-\frac{1}{q}\right)_+
\end{equation}
and argue as follows. Let $c(2)$ stand for the constant $c_2$ from Proposition \ref{pro:2} and set
\begin{equation*}
n(1)=\lceil c(2)n(0)\log(n(0)+1)\rceil.
\end{equation*}
Then for each $n\in\N$, $n\ge n(1)$ there is a unique $k(n)\ge n(0)$ so that
\begin{equation}
\label{AK8}
\lceil c(2)k(n)\log(k(n)+1)\rceil\le n< \lceil c(2)(k(n)+1)\log(k(n)+2)\rceil.
\end{equation}
Since  by \eqref{AK4} $n(0)\ge 2$, there are constants $c_{1-4}>0$ such that for all $n\in\N$, $n\ge n(1)$  
\begin{eqnarray}
c_1\log(n+1)&\le & \log(k(n)+1)\le c_2\log(n+1)
\label{AK5}\\
c_3n(\log(n+1))^{-1}&\le& k(n)\le  c_4n(\log(n+1))^{-1}.
\label{AK6}
\end{eqnarray}
By \eqref{AH0} of Proposition \ref{pro:2}
$$
\ca(A_{k(n),\omega}^{(4)},\mathcal{W}_p^r(D))\le c_2k(n)(\log(k(n)+1))
$$
and therefore monotonicity of the minimal errors, \eqref{AK8}, \eqref{AG9}, \eqref{AH9}, \eqref{AK5}, and \eqref{AK6} imply
\begin{eqnarray*}
\lefteqn{e_n^\ranno(S,B_{\mathcal{W}_p^r(D)},L_q(D_1))}
\\
&\le& e_{\left\lceil c(2)k(n)(\log(k(n)+1))^{\beta_1}\right\rceil}^\ranno(S,B_{\mathcal{W}_p^r(D)},L_q(D_1))
\le   c\Phi_1(k(n)) (\log(k(n)+1))^{\beta_1}
\\
&\le&   ck(n)^{-\frac{r}{d_1}+\left(\frac{1}{p}-\frac{1}{q}\right)_+} (\log(k(n)+1))^{(\sigma_1+1)\left(\frac{r}{d_1}-\left(\frac{1}{p}-\frac{1}{q}\right)_+\right)+1}
\\
&\le&   cn^{-\frac{r}{d_1}+\left(\frac{1}{p}-\frac{1}{q}\right)_+} (\log(n+1))^{\sigma_1\left(\frac{r}{d_1}-\left(\frac{1}{p}-\frac{1}{q}\right)_+\right)+1}
\\
&=& c\Phi_1(n)(\log(n+1))^{\frac{r}{d_1}-\left(\frac{1}{p}-\frac{1}{q}\right)_+ +1}
=c\Phi_1(n)(\log(n+1))^{2-\frac{1}{\bar{p}}}.
\end{eqnarray*}
This completes the proof of the upper bounds of \eqref{mc-eq:2c} and \eqref{B6} for $n\ge n(1)$, for $n<n(1)$ they follow trivially from the continuity of $S$.

In the adaptive case we argue similarly. Let $c(3)$ denote the constant $c_2$ from Proposition \ref{pro:3}. 
Setting
\begin{equation*}
n(2)=\lceil c(3)n(0)(\log(n(0)+1))^{\beta_2+1}\rceil,
\end{equation*}
it follows that for each $n\in\N$ with $n\ge n(2)$ there is a unique $k(n)\ge n(0)$ so that
\begin{equation}
\label{AL0}
\left\lceil c(3)k(n)(\log(k(n)+1))^{\beta_2+1}\right\rceil\le n< \left\lceil c(3)(k(n)+1)(\log(k(n)+2))^{\beta_2+1}\right\rceil.
\end{equation}
Furthermore there are constants $c_{1-4}>0$ such that  
\begin{eqnarray}
c_1\log(n+1)&\le & \log(k(n)+1)\le c_2\log(n+1)
\label{AL3}\\
c_3n(\log(n+1))^{-(\beta_2+1)}&\le& k(n)\le  c_4n(\log(n+1))^{-(\beta_2+1)}.
\label{AL2}
\end{eqnarray}
By \eqref{AH2} of Proposition \ref{pro:3}
$$
\ca(A_{k(n),\omega}^{(5)},\mathcal{W}_p^r(D))\le c_2k(n)(\log(k(n)+1))^{\beta_2+1}
$$
and therefore monotonicity of the minimal errors, \eqref{AL0} and \eqref{AH1} imply
\begin{eqnarray}
\lefteqn{e_n^\ran(S,B_{\mathcal{W}_p^r(D)},L_q(D_1))}
\notag\\
&\le& e_{\left\lceil c(3)k(n)(\log(k(n)+1))^{\beta_2+1}\right\rceil}^\ran(S,B_{\mathcal{W}_p^r(D)},L_q(D_1))
\le   c\Phi_2(k(n)) (\log(k(n)+1))^{\beta_2}.
\label{AI0}
\end{eqnarray}
If $\beta_2=0$, this together with  \eqref{AL2}    gives
\begin{equation*}
e_n^\ran(S,B_{\mathcal{W}_p^r(D)},L_q(D_1))\le  c\Phi_2\left(\frac{n}{\log(n+1)}\right).
\end{equation*}
If $\beta_2=1$, it follows from \eqref{C0} and \eqref{Y5} that
$$
\Phi_2\left(k(n)\right)= k(n)^{-\frac{r}{d_1}+\frac{1}{p}-\frac{1}{q}},
$$
which combined with \eqref{AL3}, \eqref{AL2}, and \eqref{AI0}  yields
\begin{eqnarray*}
e_n^\ran(S,B_{\mathcal{W}_p^r(D)},L_q(D_1))&\le& ck(n)^{-\frac{r}{d_1}+\frac{1}{p}-\frac{1}{q}} \log(k(n)+1)
\\
&\le&   cn^{-\frac{r}{d_1}+\frac{1}{p}-\frac{1}{q}} (\log(n+1))^{2\left(\frac{r}{d_1}-\frac{1}{p}+\frac{1}{q}\right)+1}
\\
&=&  c\Phi_2\left(\frac{n}{\log(n+1)}\right) (\log(n+1))^{\frac{r}{d_1}-\frac{1}{p}+\frac{1}{q}+1},
\end{eqnarray*}
which completes the proof of  \eqref{B9} for $n\ge n(2)$, since by \eqref{Y5}  $\beta_2=1$ implies $\frac{r}{d_1}=1-\frac{1}{q}$. The case $n<n(2)$ is trivial.

\medskip
Now we prove the lower bounds. 
For $\iota=1,2$ let $\psi^{({\iota})}$ be a $C^{\infty}$ function on $\R^{d_{\iota}}$ 
with $\text{supp}\; \psi^{({\iota})} \subset (0,1)^{d_{\iota}}$ and
$$
\tau_{\iota}:= \int_{D_{\iota}} \psi^{({\iota})}(t)\,dt\ne 0.
$$ 
For $n\ge n(3)$ we define
\begin{eqnarray}
&&l_0(n) = \left\lceil \frac{\log n-\log c(0)+1}{d_1+d_2} \right\rceil, \quad l_1(n) = \left\lceil \frac{(d_1+d_2)l_0(n)-\sigma_1\log l_0(n) }{d_1}\right\rceil\label{Y1A}
\end{eqnarray}
where $c(0)$ stands for the constant $c_0$ from Proposition \ref{pro:1}, while the constant $n(3)\in\N$ is chosen in such a way that for $n\ge n(3)$ 
\begin{equation}
\label{Z6A}
2\le l_0(n)<l_1(n).
\end{equation}
The lower bounds for $n< n(3)$ follow by monotonicity from those for $n\ge n(3)$, so in the sequel we assume $n\ge n(3)$. 
Similar to \eqref{AC7}--\eqref{AD9} we derive from \eqref{Y1A} and \eqref{Z6A}
\begin{eqnarray}
c(0)^{-\frac{1}{d_1+d_2}}n^{\frac{1}{d_1+d_2}}&<& 2^{l_0(n)} \le cn^{\frac{1}{d_1+d_2}},\label{AJ8}
\\
c_1\log(n+1)&\le& l_0(n)<l_1(n)\le c_2\log(n+1).
\label{AK3}\\
c_3n^{\frac{1}{d_1}}(\log(n+1))^{-\frac{\sigma_1}{d_1}}&\le& 2^{l_1(n)}\le c_4n^{\frac{1}{d_1}}(\log(n+1))^{-\frac{\sigma_1}{d_1}}.\label{AK2}
\end{eqnarray}
Let $l\in \{l_0(n),l_1(n)\}$. We put
\begin{equation}
\label{eq:83}
N_{1,l}=2^{d_1l},\quad N_{2,l} = 2^{d_2l}
\end{equation}
and conclude from  \eqref{AJ8} and  \eqref{Z6A} that
\begin{eqnarray}
n&<&c(0)2^{(d_1+d_2)l_0(n)}=c(0)N_{1,l_0(n)}N_{2,l_0(n)}<N_{1,l_1(n)}N_{2,l_1(n)}.
\label{AC8}
\end{eqnarray}
Furthermore, from  \eqref{AJ8}--\eqref{eq:83} we obtain
\begin{equation}
\label{AD6}
c_1n^{\frac{d_1}{d_1+d_2}}\le N_{1,l_0(n)}\le c_2n^{\frac{d_1}{d_1+d_2}},\quad c_1n(\log(n+1))^{-\sigma_1}\le N_{1,l_1(n)}\le c_2n(\log(n+1))^{-\sigma_1}.
\end{equation}
Put for $1\le i_1\le 2^{d_1l}$ and  $1\le i_2\le2^{d_2l}$
\begin{equation*}
\psi_{li_1}^{(1)} = R_{l i_1}^{(1)}\psi^{(1)},\quad \psi_{li_2}^{(2)} = R_{l i_2}^{(2)}\psi^{(2)},
\end{equation*}
with  $R_{l i}^{(1)}$ and $R_{l j}^{(2)}$  as defined before (see (\ref{WH4})),
and 
$$
\psi_{li}(s,t)= \psi_{li_1}^{(1)}(s)\psi_{li_2}^{(2)}(t)\quad (i=2^{d_2l}(i_1-1)+i_2,\,s\in [0,1]^{d_1},\,t\in [0,1]^{d_2}).
$$
We have
\begin{eqnarray}
(S\psi_{li})(s)&=&\int_{D_{l i_2}^{(2)}} \psi_{li}(s,t)\,dt 
=  \tau_2N_{2,l}^{-1}\psi_{li_1}^{(1)}(s) \quad(s\in D_1),
\label{eq:84}\\
\| \psi_{li} \|_{\mathcal{W}_p^r(D)} &\le& c2^{rl-(d_1+d_2)l/p}=c2^{rl}(N_{1,l}N_{2,l})^{-1/p}.\label{eq:85}
\end{eqnarray}
Define 
\begin{eqnarray}
&&\Gamma_l: L_p^{N_{1,l},N_{2,l}} \to \mathcal{W}_p^r(D), \quad \Gamma_l f =\sum_{i=1}^{N_{1,l}N_{2,l}}  f(i)\psi_{li}, 
\label{T9}\\
&&\Psi_l: L_q(D_1) \rightarrow L_q^{N_{1,l}},\quad (\Psi_l g)(i_1) = N_{1,l} \int_{D_{li_1}^{(1)}} g(s)ds \quad (1\le i_1\le N_{1,l}).
\label{U0}
\end{eqnarray}
Then 
\begin{equation}
\label{AD0}
\|\Gamma_l:L_p^{N_{1,l},N_{2,l}} \to \mathcal{W}_p^r(D)\|\le c2^{rl}, \quad \|\Psi_l: L_q(D_1) \rightarrow L_q^{N_{1,l}}\|\le c.
\end{equation}
Indeed, by (\ref{eq:85}), for $p<\infty$
\begin{equation*}
\|\Gamma_l f\|_{\mathcal{W}_p^r(D)}= \left( \sum_{i=1}^{N_{1,l}N_{2,l}} |f(i)|^p \| \psi_{li}\|_{\mathcal{W}_p^r(D)}^p
\right)^{1/p} 
\leq c  2^{rl}\| f\|_{L_p^{N_{1,l},N_{2,l}}},
\end{equation*}
and by H\"older's inequality, for $q<\infty$.
\begin{eqnarray*}
\| \Psi_l g \|_{L_q^{N_{1,l}}}^q & = &
N_{1,l}^{q-1} \sum_{i_1=1}^{N_{1,l}}\bigg| \int_{D_{li_1}^{(1)}} g(s)ds \bigg|^{q}\nonumber \\
& \leq &
N_{1,l}^{q-1} \sum_{i_1=1}^{N_{1,l}}\int_{D_{li_1}^{(1)}} |g(s)|^q ds \big|D_{li_1}^{(1)}\big|^{q-1} 
         = \| g\|_{L_q(D_1)}^q. 
\end{eqnarray*}
The cases $p=\infty$ and $q=\infty$ are analogous.
Moreover, from (\ref{eq:84}),  \eqref{T9}, and \eqref{U0} we conclude, identifying $i$ with $(i_1,i_2)$,
\begin{eqnarray}
\label{eq:89}
\Psi_l  S  \Gamma_l f 
& = &
\Psi_l  S \sum_{i_1=1}^{N_{1,l}}\sum_{i_2=1}^{N_{2,l}}  f(i_1,i_2)\psi_{lij}=
\Psi_l  \sum_{i=1}^{N_{1,l}}\left(\tau_2  N_{2,l}^{-1}\sum_{j=1}^{N_{2,l}} 
f(i,j) \right) \psi_{li}^{(1)} 
\notag \\
& = & \tau_1 \tau_2\sum_{i_1=1}^{N_{1,l}}\left( N_{2,l}^{-1}\sum_{i_2=1}^{N_{2,l}}
 f(i_1,i_2) \right) e_{i} = \tau_1\tau_2 S^{N_{1,l},N_{2,l}}(f),
\end{eqnarray}
with $e_i$ the unit vectors in $\K^{N_{1,l}}$.
By \eqref{T9}--\eqref{AD0},   $\Psi_l S  \Gamma_l$ reduces to $S$ and  
\begin{eqnarray}
\label{eq:90}
e_{n}^{\rm set}(\Psi_l S  \Gamma_l , B_{L_p^{N_{1,l},N_{2,l}}},L_q^{N_{1,l}})& \leq & 
e_{n}^{\rm set}(S, c 2^{rl}B_{\mathcal{W}_p^r(D)},L_q(D_1)) 
\notag\\
 &=& 
c 2^{rl}e_{n}^{\rm set}(S,B_{\mathcal{W}_p^r(D)},L_q(D_1)) ,
\end{eqnarray}
where ${\rm set}\in\{\de, \text{det--non},\ran,\text{ran--non}\}$, see \cite{Hei05b}, Proposition 1, for details on reduction. These results were shown in \cite{Hei05b} for the adaptive setting. The non-adaptive case is much easier, essentially straight-forward.
Combining \eqref{eq:89} and \eqref{eq:90}, we get
\begin{equation}
\label{eq:91}
e_n^{\rm set}(S,B_{\mathcal{W}_p^r(D)},L_q(D_1)) \geq c\,2^{-rl}
e_{n}^{\rm set}(S^{N_{1,l},N_{2,l}}, B_{L_p^{N_{1,l},N_{2,l}}},L_q^{N_{1,l}}).
\end{equation}
We conclude from \eqref{AD2} and \eqref{AD4} of Proposition \ref{pro:1} and  \eqref{AC8} that for $l\in \{l_0(n),l_1(n)\}$
\begin{eqnarray}
\label{AD1}
\lefteqn{e_n^{\rm{set}_1}(S,B_{\mathcal{W}_p^r(D)},L_q(D_1))}
\nonumber\\
&\ge&
 c\,2^{-rl}N_{1,l}^{\left(1/p-1/q\right)_+}\left\lceil\frac{n}{N_{1,l}}\right\rceil^{-\left(1-1/\bar{p}\right)}\left( \min \left( \log (N_{1,l}+1), \left\lceil\frac{n}{N_{1,l}}\right\rceil \right) \right)^{\sigma_1/2},
\end{eqnarray}
where 
\begin{equation*}
\begin{array}{lll}
  \rm{set}_1 \in \{\ran,\text{ran--non}\}& \quad\mbox{if}\quad (p\le 2)\vee (p\ge q)   \\[.2cm]
   \rm{set}_1 = \text{ran--non}& \quad\mbox{if}\quad 2<p<q.    
    \end{array}
\end{equation*}
Then \eqref{eq:83} and \eqref{AD6} imply
\begin{eqnarray}
\min \left( \log (N_{1,l_0(n)}+1), \left\lceil\frac{n}{N_{1,l_0(n)}}\right\rceil \right)
&\ge& c\min \Big(l_0(n), n^{\frac{d_2}{d_1+d_2}}\Big)\ge c\log(n+1),
\label{AD7}\\
\min \left( \log (N_{1,l_1(n)}+1), \left\lceil\frac{n}{N_{1,l_1(n)}}\right\rceil \right)
&\ge& c\min \Big(l_1(n), (\log(n+1))^{\sigma_1}\Big)
\notag\\
&\ge& c(\log(n+1))^{\sigma_1}.
\label{AD8}
\end{eqnarray}
From \eqref{AD1}, \eqref{AD6}, \eqref{AD7}, \eqref{AD8}, and \eqref{T1}, we obtain  
\begin{eqnarray*}
e_n^{\rm{set}_1}(S,B_{\mathcal{W}_p^r(D)},L_q(D_1))
&\ge& cn^{-\frac{r-\left(\frac{1}{p}-\frac{1}{q}\right)_+d_1+\left(1-\frac{1}{\bar{p}}\right)d_2}{d_1+d_2}}(\log(n+1))^{\frac{\sigma_1}{2}}
\notag\\
&&+cn^{-\frac{r}{d_1}+\left(\frac{1}{p}-\frac{1}{q}\right)_+} \,(\log(n+1))^{\sigma_1\left(\frac{r}{d_1}-\left(\frac{1}{p}-\frac{1}{q}\right)_+-\left(1-\frac{1}{\bar{p}}\right)+\frac{1}{2}\right)}\ge c\Phi_1(n),
\end{eqnarray*}
since $\sigma_1=1$ implies $\bar{p}=2$. This proves the lower bounds in \eqref{mc-eq:2c} and \eqref{B6}. 

Now we turn to the adaptive randomized setting  in the case $2<p<q$, so $\sigma_1 = 0$. Here \eqref{AD3} of Proposition \ref{pro:1}, \eqref{AC8}, and \eqref{eq:91} give for $l\in \{l_0(n),l_1(n)\}$
\begin{eqnarray*}
\lefteqn{e_n^\ran(S,B_{\mathcal{W}_p^r(D)},L_q(D_1))}
\nonumber\\ 
&\ge&c2^{-rl}N_{1,l}^{1/p-1/q}\left\lceil\frac{n}{N_{1,l}}\right\rceil^{-\left(1-1/p\right)}+c2^{-rl}\left\lceil\frac{n}{N_{1,l}}\right\rceil^{-1/2}  \left( \log (N_{1,l}+1) \right)^{\delta_{q,\infty}/2}
\end{eqnarray*}
and therefore, using  \eqref{AD6},
\begin{eqnarray}
\lefteqn{e_n^\ran(S,B_{\mathcal{W}_p^r(D)},L_q(D_1))}
\notag\\
&\ge&c\,2^{-rl_0(n)}N_{1,l_0(n)}^{1/p-1/q}\left\lceil\frac{n}{N_{1,l_0(n)}}\right\rceil^{-\left(1-1/p\right)}+c\,2^{-rl_0(n)}\left\lceil\frac{n}{N_{1,l_0(n)}}\right\rceil^{-1/2}  \left( \log (N_{1,l_0(n)}+1) \right)^{\delta_{q,\infty}/2}
\notag\\
&&+c\,2^{-rl_1(n)}N_{1,l_1(n)}^{1/p-1/q}\left\lceil\frac{n}{N_{1,l_1(n)}}\right\rceil^{-\left(1-1/p\right)}+c\,2^{-rl_1(n)}\left\lceil\frac{n}{N_{1,l_1(n)}}\right\rceil^{-1/2}  \left( \log (N_{1,l_1(n)}+1) \right)^{\delta_{q,\infty}/2}
\notag\\
&\ge&cn^{\frac{-r+\big(\frac{1}{p}-\frac{1}{q}\big)d_1-\big(1-\frac{1}{p}\big)d_2 }{d_1+d_2}}
+cn^{\frac{-r-\frac12 d_2}{d_1+d_2}}(\log(n+1))^{\delta_{q,\infty}/2}
+cn^{-\frac{r}{d_1}+\frac{1}{p}-\frac{1}{q}}\ge c\Phi_2(n),
\label{AG8}
\end{eqnarray}
which is the lower bound of \eqref{B9} for $\sigma_2=0$. Now assume $\sigma_2=1$, which by \eqref{AE1} implies $q=\infty$ and
\begin{equation}
\left(\bigg(\frac{1}{2}-\frac{1}{p}\bigg)d_2\ge \bigg(\frac{1}{p}-\frac{1}{q}\bigg)d_1\right)\wedge \left(\frac{r}{d_1}\ge  \left(\frac{1}{p}-\frac{1}{q}\right)\left(\frac{d_1}{d_2}+1\right)+\frac{1}{2}\right).
\label{AF3}
\end{equation}
As mentioned after relation \eqref{C6A}, relations \eqref{B7} and \eqref{C7} also hold with ''$>$'' replaced by ''$\ge$''  
so we conclude from them and \eqref{AF3} that $\Phi_2(n)=n^{\frac{-r-\frac12 d_2}{d_1+d_2}}$ and therefore from \eqref{AG8}
\begin{eqnarray*}
e_n^\ran(S,B_{\mathcal{W}_p^r(D)},L_q(D_1))
&\ge&
cn^{\frac{-r-\frac12 d_2}{d_1+d_2}}(\log(n+1))^{1/2}\ge c\Phi_2(n)(\log(n+1))^{1/2}, 
\end{eqnarray*}
showing \eqref{B9} also for the case $\sigma_2=1$.

In the deterministic case we have by \eqref{eq:91}, \eqref{AD5}, \eqref{AJ8}, and \eqref{AD6}
\begin{equation*}
e_n^{\de}(S,B_{\mathcal{W}_p^r(D)},L_q(D_1))\ge
 c\,2^{-rl_0(n)}N_{1,l_0(n)}^{\left(1/p-1/q\right)_+}\ge cn^{\frac{-r+d_1\left(\frac{1}{p}-\frac{1}{q}\right)_+}{d_1+d_2}},
\end{equation*}
which is the lower bound of \eqref{AE2}.
This completes the proof of the lower bounds and thus of the theorem.

\end{proof}

\section{Speedup}
\label{sec:6}
In this section we assume that $r\in\N_0$, $d_1,d_2\in\N$, $2<p<q\le \infty$ and \eqref{B1} is satisfied.
As in Part I we want to determine  the widest  gap between non-adaptive and adaptive randomized minimal errors. We define the gap as
\begin{equation*}
\gamma(n;r,p,q,d_1,d_2)= \frac{e_n^\ranno(S, B_{\mathcal{W}_p^r(D)},L_q(D_1))}{e_n^\ran(S, B_{\mathcal{W}_p^r(D)},L_q(D_1))}
\end{equation*}
and the principal exponent of  $\gamma(n;r,p,q,d_1,d_2)$: 
\begin{equation}
\label{AJ0}
\theta(r,p,q,d_1,d_2)= \lim_{n\to \infty}\frac{\log\gamma(n;r,p,q,d_1,d_2)}{\log n}.
\end{equation}
\begin{corollary}
\label{cor:3}
Under the assumptions above the limit in \eqref{AJ0} exists and
\begin{equation*}
\max\theta(r,p,q,d_1,d_2)=\frac18,
\end{equation*}
where the maximum is taken over all $r,d_1,d_2,p,q$ fulfilling the assumption.
The maximum is attained iff 
\begin{equation}
\label{AJ3}
p=4,\quad q=\infty,\quad r\ge d_1=d_2.
\end{equation}
In this case there are constants $c_1,c_2>0$ such that for all $n\in\N$
\begin{eqnarray}
c_1n^{\frac{1}{8}}(\log(n+1))^{-\frac{5}{2}}
&\le& \gamma\left(n;r,4,\infty,d_1,d_1\right)\le c_2n^{\frac{1}{8}}(\log(n+1))^{-\frac{1}{2}}\quad \text{if}\quad r=d_1
\label{AJ5}\\
c_1n^{\frac{1}{8}}(\log(n+1))^{-\frac{r}{2d_1}-\frac{1}{4}}
&\le& \gamma\left(n;r,4,\infty,d_1,d_1\right)\le c_2n^{\frac{1}{8}}(\log(n+1))^{-\frac{1}{2}}\quad \text{if}\quad r>d_1.
\label{AJ6}
\end{eqnarray}
\end{corollary}
\begin{proof} By Theorem \ref{theo:1} we have 
\begin{equation}
c_1\frac{\Phi_1(n)}  {\Phi_2\left(\frac{n}{\log(n+1)}\right)(\log(n+1))^{\beta_2\left(2-\frac{1}{p}\right)}}
\le \gamma(n)\le c_2\frac{\Phi_1(n)(\log (n+1))^{\beta_1\left(2-\frac{1}{p}\right)}}{\Phi_2(n)(\log(n+1))^{\sigma_2/2}},
\label{AJ1}
\end{equation}
which together with the definitions \eqref{T1}, \eqref{C7}, and  \eqref{C0} of $\Phi_1$  and $\Phi_2$ shows that the limit in \eqref{AJ0} exists. Now we determine $\theta$. 
First we  assume \eqref{B7}, 
which is equivalent to 
\begin{equation}
\frac{d_1+d_2}{d_1} > \frac{\frac{1}{2}-\frac{1}{p}+\frac{1}{p}-\frac{1}{q}}{\frac{1}{2}-\frac{1}{p}}=\frac{\frac{1}{2}-\frac{1}{q}}{\frac{1}{2}-\frac{1}{p}}.
\label{AI6}
\end{equation}
Then \eqref{T1}, \eqref{C7}, and \eqref{AJ1} imply the following.
\begin{eqnarray*}
&&\text{ If } \quad \frac{r}{d_1}\le \frac{1}{p}-\frac{1}{q}+\frac{1}{2},\quad \text{then}\quad\theta=0, 
\notag\\[.4cm]
&&\text{ if }\quad\frac{1}{p}-\frac{1}{q}+\frac{1}{2}< \frac{r}{d_1}\le \left(\frac{1}{p}-\frac{1}{q}\right)\left(\frac{d_1}{d_2}+1\right)+\frac{1}{2},\quad \text{then}
\notag\\[.2cm]
&&\qquad\quad\theta=\frac{-r+\left(\frac{1}{p}-\frac{1}{q}\right)d_1-\frac{d_2}{2}} {d_1+d_2}-\left(-\frac{r}{d_1}+\frac{1}{p}-\frac{1}{q}\right)=\frac{\left(\frac{r}{d_1}-\left(\frac{1}{p}-\frac{1}{q}+\frac{1}{2}\right)\right)d_2} {d_1+d_2},\quad
\\[.4cm]
&&\text{ if } \quad    \left(\frac{1}{p}-\frac{1}{q}\right)\left(\frac{d_1}{d_2}+1\right)+\frac{1}{2}<\frac{r}{d_1},\quad \text{then}
\notag\\[.2cm]
&&\qquad\quad\theta=\frac{-r+\left(\frac{1}{p}-\frac{1}{q}\right)d_1-\frac{d_2}{2}} {d_1+d_2}-\frac{-r-\frac12 d_2}{d_1+d_2}
 =   \frac{\left(\frac{1}{p}-\frac{1}{q}\right)d_1} {d_1+d_2}.
\notag
\end{eqnarray*}
Clearly,
$$
\frac{r}{d_1}\le \left(\frac{1}{p}-\frac{1}{q}\right)\left(\frac{d_1}{d_2}+1\right)+\frac{1}{2}
$$
implies
$$
\frac{\left(\frac{r}{d_1}-\left(\frac{1}{p}-\frac{1}{q}+\frac{1}{2}\right)\right)d_2} {d_1+d_2}
\le\frac{\left(\frac{1}{p}-\frac{1}{q}\right)d_1} {d_1+d_2}.
$$
Consequently, taking into account \eqref{AI6}, we obtain
\begin{equation}
\label{AI8}
\theta\le\frac{\left(\frac{1}{p}-\frac{1}{q}\right)d_1} {d_1+d_2}<\frac{\left(\frac{1}{p}-\frac{1}{q}\right)\left(\frac{1}{2}-\frac{1}{p}\right)}{\frac{1}{2}-\frac{1}{q}}\le \frac{\frac{1}{2}-\frac{1}{q}}{4}\le \frac{1}{8}.
\end{equation}

Now we  assume \eqref{B8}, 
which is equivalent to 
\begin{equation*}
\frac{d_1+d_2}{d_2} \ge \frac{\frac{1}{2}-\frac{1}{q}}{\frac{1}{p}-\frac{1}{q}}.
\end{equation*}
By \eqref{T1}, \eqref{C0}, and \eqref{AJ1} the following is true.
\begin{eqnarray*}
&&\text{ If } \quad \frac{r}{d_1}\le \frac{1}{p}-\frac{1}{q}+\frac{1}{2},\quad \text{then}\quad\theta=0, 
\notag\\[.4cm]
&&\text{ if }\quad\frac{1}{p}-\frac{1}{q}+\frac{1}{2}< \frac{r}{d_1}\le 1-\frac{1}{q},\quad \text{then}
\\[.2cm]
&&\qquad\quad\theta=\frac{-r+\left(\frac{1}{p}-\frac{1}{q}\right)d_1-\frac{d_2}{2}} {d_1+d_2}-\left(-\frac{r}{d_1}+\frac{1}{p}-\frac{1}{q}\right)=\frac{\left(\frac{r}{d_1}-\left(\frac{1}{p}-\frac{1}{q}+\frac{1}{2}\right)\right)d_2} {d_1+d_2},\quad
\\[.4cm]
&&\text{ if } \quad    1-\frac{1}{q}<\frac{r}{d_1},\quad \text{then}
\notag\\[.2cm]
&&\qquad\quad\theta=\frac{-r+\left(\frac{1}{p}-\frac{1}{q}\right)d_1-\frac{d_2}{2}} {d_1+d_2}-\frac{-r+\big(\frac{1}{p}-\frac{1}{q}\big)d_1-\big(1-\frac{1}{p}\big)d_2 }{d_1+d_2}=\frac{\left(\frac{1}{2}-\frac{1}{p}\right)d_2} {d_1+d_2}.
\notag
\end{eqnarray*}
For $\frac{1}{p}-\frac{1}{q}+\frac{1}{2}<\frac{r}{d_1}< 1-\frac{1}{q}$ we have 
$$
\theta=\frac{\left(\frac{r}{d_1}-\left(\frac{1}{p}-\frac{1}{q}+\frac{1}{2}\right)\right)d_2} {d_1+d_2}
<\frac{\left(\frac{1}{2}-\frac{1}{p}\right)d_2} {d_1+d_2},
$$
while for $\frac{r}{d_1}\ge 1-\frac{1}{q}$ 
$$
\theta= \frac{\left(\frac{1}{2}-\frac{1}{p}\right)d_2} {d_1+d_2}\le\frac{\left(\frac{1}{2}-\frac{1}{p}\right)\left(\frac{1}{p}-\frac{1}{q}\right)}{\frac{1}{2}-\frac{1}{q}}\le \frac{1}{8},
$$
see \eqref{AI8}. We conclude that the maximal value of $\theta$ is reached iff $\theta=\frac{1}{8}$ iff \eqref{AJ3} holds.

So assume \eqref{AJ3}. Here we estimate the gap including log factors. By \eqref{Q8}, \eqref{S9}, \eqref{Y5}, and \eqref{AE1} we have
\begin{equation*}
\sigma_1=0,\quad\beta_1=0,\quad\beta_2=\left\{\begin{array}{lll}
 1  & \quad\mbox{if}\quad  r=d_1,  \\
 0  & \quad\mbox{if}\quad   r>d_1,
    \end{array}
\right.  \quad\sigma_2=1.
\end{equation*}
From \eqref{T1} and \eqref{C0} we conclude 
\begin{equation*}
\Phi_1(n)=n^{-\frac{r}{2d_1}-\frac{1}{8}},\quad\Phi_2(n)=n^{-\frac{r}{2d_1}-\frac{1}{4}},
\end{equation*}
and consequently
\begin{equation*}
\frac{\Phi_1(n)}{\Phi_2\left(\frac{n}{\log(n+1)}\right)(\log(n+1))^{\beta_2\left(2-\frac{1}{p}\right)}}
=n^{\frac{1}{8}}(\log(n+1))^{-\frac{7\beta_2}{4}-\frac{r}{2d_1}-\frac{1}{4}}
\end{equation*}
\begin{equation*}
\frac{\Phi_1(n)(\log (n+1))^{\beta_1\left(2-\frac{1}{p}\right)}}{\Phi_2(n)(\log(n+1))^{\sigma_2/2}}
=n^{\frac{1}{8}}(\log(n+1))^{-\frac{1}{2}}.
\end{equation*}
This together with \eqref{AJ1} yields \eqref{AJ5} and  \eqref{AJ6}.

\end{proof}

\end{document}